\documentclass[11pt]{article}
\usepackage{graphicx}
\usepackage{amsmath,amsfonts,amssymb}
\usepackage{citesort}
\usepackage{url}
\usepackage{amsmath}
\usepackage{graphicx}
\usepackage{epsfig}
\usepackage{epstopdf}
\usepackage{color}
\def\disp{\displaystyle}

\def\tto{\;{\lower 1pt \hbox{$\rightarrow$}}\kern -10pt
\hbox{\raise 2pt \hbox{$\rightarrow$}}\;}

\def\Bar{\overline}
\def\ra{\rangle}
\def\la{\langle}

\def\B{\Bbb B}
\def\h{\hfill\Box}
\def\R{\Bbb R}
\def\N{\Bbb N}
\def\ox{\bar{x}}

\def\oy{\bar{y}}
\def\oz{\bar{z}}

\def\ou{\bar{u}}

\def\co{\mbox{\rm co}}
\def\ri{\mbox{\rm ri}}
\def\gph{\mbox{\rm gph}}
\def\aff{\mbox{\rm aff}}
\def\epi{\mbox{\rm epi}}

\def\dom{\mbox{\rm dom}}

\def\h{\hfill\square}

\def\O{\Omega}
\def\ph{\varphi}
\def\emp{\emptyset}

\def\lm{\lambda}

\def\dd{\delta}
\def\al{\alpha}
\def\Th{\Theta}

\newcounter{lk}

\begin{document}
\begin{center}
{\bf GEOMETRIC APPROACH TO \\CONVEX SUBDIFFERENTIAL CALCULUS}\\[1ex]
\today\\[2ex]
BORIS S. MORDUKHOVICH\footnote{Department of Mathematics, Wayne State University, Detroit, MI 48202, USA(boris@math.wayne.edu).
Research of this author was partly supported by the National Science Foundation under grants DMS-1007132 and DMS-1512846 and the Air Force Office of Scientific Research grant \#15RT0462.} and NGUYEN MAU NAM\footnote{Fariborz Maseeh Department of Mathematics and Statistics, Portland State University, PO Box 751, Portland, OR 97207, USA(mau.nam.nguyen@pdx.edu). The research of this author was partially supported by the NSF under grant DMS-1411817 and the Simons Foundation under grant \#208785.}\\[2ex]
{\bf Dedicated to Franco Giannessi and Diethard Pallaschke with great respect}
\end{center}
\small{\bf Abstract.} In this paper we develop a geometric approach to convex subdifferential calculus in finite dimensions with employing some ideas of modern variational analysis. This approach allows us to obtain natural and rather easy proofs of basic results of convex subdifferential calculus in full generality and also derive new results of convex analysis concerning optimal value/marginal functions, normals to inverse images of sets under set-valued mappings, calculus rules for coderivatives of single-valued and set-valued mappings, and calculating coderivatives of solution maps to parameterized generalized equations governed by set-valued mappings with convex graphs.

\newtheorem{Theorem}{Theorem}[section]
\newtheorem{Proposition}[Theorem]{Proposition}
\newtheorem{Remark}[Theorem]{Remark}
\newtheorem{Lemma}[Theorem]{Lemma}
\newtheorem{Corollary}[Theorem]{Corollary}
\newtheorem{Definition}[Theorem]{Definition}
\newtheorem{Example}[Theorem]{Example}
\renewcommand{\theequation}{\thesection.\arabic{equation}}
\normalsize

\section{Introduction}
\setcounter{equation}{0}

The notion of {\em subdifferential} (collection of {\em subgradients}) for nondifferentiable convex functions was independently introduced and developed by Moreau \cite{mor} and Rockafellar \cite{roc} who were both influenced by Fenchel \cite{f}. Since then this notion has become one of the most central concepts of convex analysis and its various applications, including first of all convex optimization. The underlying difference between the standard derivative of a differentiable function and the subdifferential of a convex function at a given point is that the subdifferential is a {\em set} (of subgradients) which reduces to a singleton (gradient) if the function is differentiable. Due to the set-valuedness of the subdifferential, deriving {\em calculus rules} for it is a significantly more involved task in comparison with classical differential calculus. Needless to say that subdifferential calculus for convex functions is at the same level of importance as classical differential calculus, and it is difficult to imagine any usefulness of subgradients unless reasonable calculus rules are available.

The first and the most important result of convex subdifferential calculus is the subdifferential {\em sum rule}, which was obtained at the very beginning of convex analysis and has been since known as the {\em Moreau-Rockafellar theorem}. The reader can find this theorem and other results of convex subdifferential calculus in finite-dimensional spaces in the now classical monograph by Rockafellar \cite{r}. More results in this direction in finite and infinite dimensions with various applications to convex optimization, optimal control, numerical analysis, approximation theory, etc. are presented, e.g., in the monographs \cite{Bauschke2011,Bertsekas2003,Borwein2000,Boyd2004,g,HU,KK2013,Tikhomirov2003,bmn,pr,p,z} among the vast bibliography on the subject. In the recent time, convex analysis  has become more and more important for applications to many new fields  such as computational statistics, machine learning, and sparse optimization. Having this in mind, our major goal here is to revisit the convex subdifferential and provide an easy way to excess basic subdifferential calculus rules in finite dimensions.

In this paper, which can be considered as a supplement to our recent book \cite{bmn}, we develop a {\em a geometric approach} to convex subdifferential calculus. Our approach relies on the normal cone {\em intersection rule} based on {\em convex separation} and derives from it the major rules of subdifferential calculus without any appeal to duality, directional derivatives, and other tangentially generated constructions. This approach allows us to give direct and simple proofs of known results of convex subdifferential calculus in full generality and also to obtain some new results in this direction as those presented in Sections~9, 11, and 12. The developed approach is largely induced by the {\em dual-space} geometric approach of general variational analysis based on the {\em extremal principle} for systems of sets, which can be viewed as a variational counterpart of convex separation without necessarily imposing convexity; see \cite{m-book1} and the references therein.

Some of the results and proofs presented below have been outlined in the exercises of our book \cite{bmn} while some other results (e.g., those related to subgradients of the optimal value function, coderivatives and their applications, etc.) seem to be new in the convex settings under consideration. In order to make the paper self-contained for the reader's convenience and also to make this material to be suitable for teaching, we recall here some basic definitions and properties of convex sets and functions with illustrative figures and examples. Our notation follows \cite{bmn}.

\section{Basic Properties of Convex Sets}
\setcounter{equation}{0}

Here we recall some basic concepts and properties of convex sets. The detailed proofs of all the results presented in this and the next section can be found in \cite{bmn}. Throughout the paper, consider the Euclidean space $\R^n$ of $n-$tuples of real numbers with the inner product
$$
\la x,y\ra:=\sum_{i=1}^n x_iy_i\;\mbox{ for }\;x=(x_1,\ldots,x_n)\in\R^n\;\mbox {and }\;y=(y_1,\ldots, y_n)\in\R^n.
$$
The Euclidean norm induced by this inner product is defined as usual by
$$
\|x\|:=\sqrt{\sum_{i=1}^n x_i^2}.
$$
We often identify each element $x=(x_1,\ldots, x_n)\in\R^n$ with the column $x=[x_1,\ldots, x_n]^\top$.

Given two points $a,b\in\R^n$, the \emph{line segment/interval} connecting $a$ and $b$ is
\begin{equation*}
[a,b]:=\{\lambda a+(1-\lambda)b\;|\;\lambda\in[0,1]\}.
\end{equation*}

A subset $\Omega$ of $\R^n$ is called \emph{convex} if $\lambda x+(1-\lambda)y\in \O$ for all $x,y\in\O$ and $\lambda\in(0,1)$. A mapping $B\colon\R^n\to\R^p$ is \emph{affine} if there exist a $p\times n$ matrix $A$ and a vector $b\in\R^p$ such that
$$
B(x)=Ax+b\;\mbox{\rm for all }\;x\in\R^n.
$$
It is easy to check that the convexity of sets is preserved under images of affine mappings.

\begin{Proposition} Let $B\colon\R^n\to\R^p$ be an affine mapping. The following properties hold:\\[1ex]
{\rm\bf (i)} If $\O$ is a convex subset of $\R^n$, then the direct image $B(\Omega)$ is a convex subset of $\R^p$.\\
{\rm\bf (ii)} If $\Theta$ is a convex subset of $\R^p$, then the inverse image $B^{-1}(\Theta)$ is a convex subset of $\R^n$.
\end{Proposition}

For any collection of convex sets $\{\O_i\}_{i\in I}$, their intersection $\bigcap_{i\in I}\O_i$ is also convex. This motivates us to define the {\em convex hull} of a given set $\O\subset\R^n$ by
\begin{equation*}
\co(\Omega):=\bigcap\Big\{C\;\Big|\;C \text{ is convex and } \Omega\subset C\Big\},
\end{equation*}
i.e., the convex hull of a set $\O$ is the smallest convex set containing $\O$. The following useful observation is a direct consequence of the definition.

\begin{Proposition}\label{cvh} For any subset $\O$ of $\R^n$, its convex hull admits the representation
\begin{equation*}
\co(\O)=\Big\{\sum_{i=1}^m\lambda_iw_i\;\Big|\;\sum_{i=1}^m\lambda_i=1,\;\lambda_i\ge 0,\;w_i\in\O,\;m\in{\Bbb N}\Big\},
\end{equation*}
where the symbol $\N$ stands for the set of all positive integers.
\end{Proposition}

Given two points $a,b\in\R^n$, the \emph{line} connecting $a$ and $b$ in $\R^n$ is defined by
\begin{equation*}
L[a,b]:=\{\lambda a+(1-\lambda)b\;|\;\lambda\in\R\}.
\end{equation*}
A subset $A$ of $\R^n$ is called \emph{affine} if for any $x,y\in A$ and for any $\lambda\in\R$ we have
$$
\lambda x+(1-\lambda)y\in A,
$$
which means that $A$ is affine if and only if the line connecting any two points $a,b\in A$ is a subset of $A$. This shows that the intersection of any collection of affine sets is an affine set and thus allows us to define the \emph{affine hull} of $\Omega$ by
\begin{equation*}
\aff(\Omega):=\bigcap\Big\{A\;\Big|\;A\;\text{ is affine and }\;\Omega\subset A\Big\}.
\end{equation*}

\begin{figure}[!ht]
\begin{center}
\includegraphics[width=8cm]{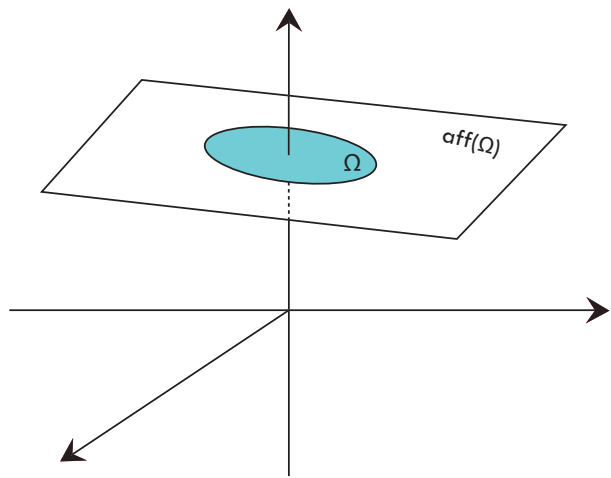}
\caption{Affine hull.}
\label{fig:1}
\end{center}
\end{figure}

Similarly to the case of the convex hull, we have the following representation.

\begin{Proposition}\label{cvh-aff} For any subset $\O$ of $\R^n$, its affine hull is represented by
\begin{equation*}
\aff(\Omega)=\Big\{\sum_{i=1}^{m}\lambda_i\omega_i\;\Big|\;\sum_{i=1}^{m}\lambda_i=1\;\omega_i \in\Omega,\;m\in\N\Big\}.
\end{equation*}
\end{Proposition}

Now we present some simple facts about affine sets.

\begin{Proposition}\label{closed property} Let $A$ be an affine subset of $\R^n$. The following properties hold:\\[1ex]
{\rm\bf (i)} If $A$ contains $0$, then it is a subspace of $\R^n$.\\
{\rm\bf (ii)} $A$ is a closed, and so the affine hull of an arbitrary set $\O$ is always closed.
\end{Proposition}
{\bf Proof.} {\bf (i)} Since $A$ is affine and since $0\in A$, for any $x\in A$ and $\lambda\in\R$ we have that $\lambda x=\lambda x+(1-\lambda)0\in A$. It also holds
\begin{equation*}
x+y=2(x/2+y/2)\in A
\end{equation*}
for any two elements $x,y\in A$, and thus $A$ is a subspace of $\R^n$.\\[1ex]
{\bf (ii)} The conclusion is obvious if $A=\emptyset$. Suppose that $A\ne\emp$, choose $x_0\in A$, and consider the set $L:=A-x_0$. Then $L$ is affine with $0\in L$, and so it is a subspace of $\R^n$. Since $A=x_0+L$ and any subspace of $\R^n$ is known to be closed, the set $A$ is closed as well. $\h$

We are now ready to formulate the notion of the \emph{relative interior} $\ri(\Omega)$ of a convex set $\O\subset\R^n$, which plays a central role in developing convex subdifferential calculus.
\begin{Definition} We say that $x\in\ri(\Omega)$ if there exists $\gamma>0$ such that
\begin{equation*}
\B(x;\gamma)\cap\aff(\Omega)\subset\Omega,
\end{equation*}
where $\B(x;\gamma)$ denotes the closed ball centered at $x$ with radius $\gamma$.
\end{Definition}

\begin{figure}[!ht]
\begin{center}
\includegraphics[width=9cm]{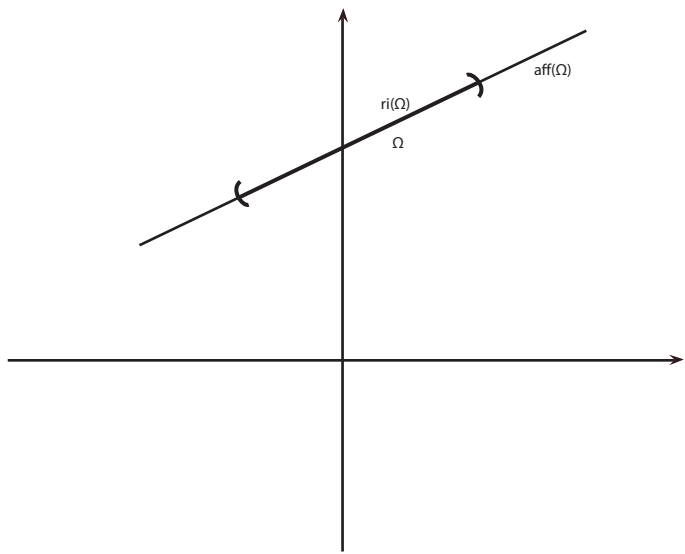}
\caption{Relative interior.}
\label{fig:1-ri}
\end{center}
\end{figure}

The following simple proposition is useful in what follows and serves as an example for better understanding of the relative interior.

\begin{Proposition}\label{prop 1.5} Let $\O$ be a nonempty convex set. Suppose that $\ox\in\ri(\O)$ and $\oy\in\O$. Then there exists $t>0$ such that
$$
\ox+t(\ox-\oy)\in\O.
$$
\end{Proposition}
{\bf Proof.} Choose a number $\gamma>0$ such that
$$
\B(\ox;\gamma)\cap\aff(\Omega)\subset\Omega
$$
and note that $\ox+t(\ox-\oy)=(1+t)\ox+(-t)\oy\in\aff(\Omega)$ for all $t\in\R$ as it is an affine combination of $\ox$ and $\oy$. Select $t>0$ so small that
$\ox+t(\ox-\oy)\in\B(\ox;\gamma)$. Then we have $\ox+t(\ox-\oy)\in\B(x;\gamma)\cap\aff(\Omega)\subset\Omega$. $\h$

Given two elements $a,b\in\R^n$, define the half-open interval
\begin{equation*}
[a,b):=\{ta+(1-t)b\;|\;0<t\le 1\}.
\end{equation*}
The following theorem is taken from \cite[Theorems~6.1 and 6.2]{r}; see also \cite[Theorem~1.72]{bmn} for a direct and detailed proof.

\begin{Theorem}\label{ri 1.72}
Let $\Omega$ be a nonempty convex subset of $\R ^n$. Then:\\
{\rm\bf (i)} We have $\ri(\Omega)\ne\emp$ and\\
{\rm\bf (ii)} $[a,b)\subset\ri(\Omega)$ for any $a\in\ri(\Omega)$ and $b\in\Bar{\Omega}$.
\end{Theorem}

The next theorem taken from \cite[Theorem~6.7]{r} gives us a convenient way to represent the relative interior of the direct image of a convex set under an affine mapping via of the relative interior of this set; see also \cite[Excercise~1.27]{bmn} and its solution for a simple proof.

\begin{Theorem}\label{ari} Let $B:\R^n\to\R^p$ be affine, and let $\O\subset\R^n$ be convex. Then we have
$$
B(\ri(\O))=\ri(B(\O)).
$$
\end{Theorem}

A useful consequence of this result is the following property concerning the {\em difference} of two subsets $A_1,A_2\subset\R^n$ defined by
$$
A_1-A_2:=\{a_1-a_2\;|\;a_1\in A_1\;\mbox{\rm and }\;a_2\in A_2\}.
$$
\begin{Corollary}\label{ri of set difference} Let $\O_1$ and $\O_2$ be convex subsets of $\R^n$. Then
$$
\ri(\O_1-\O_2)=\ri(\O_1)-\ri(\O_2).
$$
\end{Corollary}
{\bf Proof.} Consider the linear mapping $B\colon\R^n\times\R^n\to\R^n$ given by $B(x,y):=x-y$ and form the Cartesian product $\O:=\O_1\times\O_2$. Then we have
$B(\O)=\O_1-\O_2$, which yields
\begin{equation*}
\ri(\O_1-\O_2)=\ri(B(\O))=B(\ri(\O))=B(\mbox{\rm ri}(\O_1\times\O_2))=B(\ri(\O_1)\times\ri(\O_2))=\ri(\O_1)-\ri(\O_2)
\end{equation*}
by using the simple fact that $\mbox{\rm ri}(\O_1\times\O_2)=\ri(\O_1)\times\ri(\O_2)$. $\h$

Given now a set $\O\subset\Bbb R^n$, the {\em distance function} associated with $\Omega$ is defined on $\R^n$ by
\begin{equation*}
d(x;\Omega):=\inf\{\|x-\omega\|\;|\;\omega\in\Omega\}
\end{equation*}
and the {\em Euclidean projection} of $x$ onto $\O$ is
\begin{equation*}\label{e:projection}
\pi(x;\O):=\{\omega\in\O\;|\;\|x-\omega\|=d(x;\O)\}.
\end{equation*}
It is well known (see, e.g., \cite[Corollary~1.76]{bmn}) that $\pi(\ox;\O)$ is a singleton whenever the set $\O$ is nonempty, closed, and convex in $\R^n$.

The next proposition plays a crucial role in proving major results on convex separation.

\begin{Proposition}\label{angle prop} Let $\O$ be a nonempty closed convex subset of $\R^n$ with $\ox\notin\O$. Then we have $\bar\omega=\pi(\bar x;\O)$ if and only if $\bar\omega\in\O$ and
\begin{equation}\label{af}
\la\bar x-\bar\omega,\omega-\bar\omega\ra\le 0\;\mbox{ for all }\;\omega\in\O.
\end{equation}
\end{Proposition}
{\bf Proof.} Let us first show that \eqref{af} holds for $\bar\omega:=\pi(\bar x;\O)$. Fixing any $t\in(0,1)$ and $\omega\in\O$, we get $t\omega+(1-t)\bar\omega\in\O$, which implies by the projection definition that
\begin{equation*}
\|\bar x-\bar\omega\|^2\le\|\bar x-[t\omega+(1-t)\bar\omega]\|^2.
\end{equation*}
This tells us by the construction of the Euclidean norm that
\begin{equation*}
\|\bar x-\bar\omega\|^2\le\|\bar x-[\bar\omega+t(\omega-\bar\omega)\|^2=\|\bar x-\bar\omega\|^2-2t\la\bar x-\bar\omega,\omega-\bar\omega\ra+t^2\|\omega-\bar\omega\|^2
\end{equation*}
and yields therefore the inequality
\begin{equation*}
2\la\bar x-\bar\omega,\omega-\bar\omega\ra \leq t\|\omega-\bar\omega\|^2.
\end{equation*}
Letting there $t\to 0^+$ justifies property \eqref{af}.

Conversely, suppose that \eqref{af} is satisfied for $\bar\omega\in\O$ and get for any $\omega\in\O$ that
\begin{equation*}
\|\bar x-\omega\|^2=\|\bar x-\bar\omega+\bar\omega-\omega\|^2=\|\bar x-\bar\omega\|^2+2\la\bar x-\bar\omega,\bar\omega-\omega\ra+\|\bar\omega-\omega\|^2
\ge\|\bar x-\bar\omega\|^2.
\end{equation*}
Thus we have $\|\bar x-\bar\omega\|\le\|\bar x-\omega\|$ for all $\omega\in\O$, which verifies $\bar\omega=\pi(\bar x;\O)$. $\h$

\section{Basic Properties of Convex Functions}
\setcounter{equation}{0}

In this section we deal with {\em extended-real-valued} functions $f\colon\R^n\to(-\infty,\infty]=\R\cup\{\infty\}$ and use the following arithmetic conventions
on $(-\infty,\infty]$:
\begin{eqnarray*}
\begin{array}{ll}
&\alpha+\infty=\infty+\alpha=\infty\;\mbox{ for all }\;\alpha\in\R,\\
&\alpha\cdot\infty=\infty\cdot\alpha=\infty\;\mbox{\rm for all }\;\alpha>0,\\
&\infty+\infty=\infty,\quad\infty\cdot\infty=\infty,\quad 0\cdot\infty=\infty\cdot 0=0.
\end{array}
\end{eqnarray*}
The {\em domain} and {\em epigraph} of $f\colon\R^n\to(-\infty,\infty]$ are defined, respectively, by
\begin{eqnarray*}
\dom(f):=\{x\in\R^n\;|\;f(x)<\infty\},\quad\epi(f):=\big\{(x,\alpha)\in\R^{n+1}\;\big|\;x\in\R^n,\;\alpha\ge f(x)\big\}.
\end{eqnarray*}

Recall that a function $f\colon\R^n\to(-\infty,\infty]$ is {\em convex} on $\R^n$ if
\begin{equation*}
f\big(\lm x+(1-\lm)y\big)\le\lm f(x)+(1-\lm)f(y)\;\mbox{ for all }\;x,y\in\R^n\;\mbox{ and }\;\lm\in (0,1).
\end{equation*}
It is not hard to check that $f$ is convex on $\R^n$ if and only if its epigraph is a convex subset of $\R^{n+1}$. Furthermore, the domain of a convex function is a convex set.

The class of convex functions is favorable for optimization theory and applications. The next proposition reveals a characteristic feature of convex functions from the viewpoint of minimization. Recall that $f$ has a {\em local minimum} at $\ox\in\dom(f)$ if there is $\gamma>0$ such that
$$
f(x)\ge f(\ox)\;\mbox{\rm for all }\;x\in\B(\ox;\gamma).
$$
If this inequality holds for all $x\in\R^n$, we say that $f$ has an {\em absolute/global minimum} at $\ox$.

\begin{Proposition}\label{global} Let $f\colon\R^n\to(-\infty,\infty]$ be a convex function. Then $f$ has a local minimum at $\bar x$ if and only if $f$ has an absolute minimum at $\bar x$.
\end{Proposition}
{\bf Proof.} We only need to show that any local minimizer of $f$ provides a global minimum to this function. Suppose that $\ox$ is such a local minimizer, fix any $x\in\R^n$, and denote $x_k:=(1-k^{-1})\bar x+k^{-1}x$ for all $k\in\N$. Then $x_k\to\bar x$ as $k\to\infty$. Taking $\gamma>0$ from the definition of $\ox$ gives us that $x_k\in\B(\bar x;\gamma)$ when $k$ is sufficiently large. Hence
\begin{equation*}
f(\bar x)\le f(x_k)\le(1-k^{-1})f(\bar x)+k^{-1}f(x),
\end{equation*}
which readily implies that $f(\bar x)\le f(x)$ whenever $x\in\R^n$. $\h$

Next we present the basic definition of the {\em subdifferential} as the collection of {\em subgradients} for a convex function at a given point.

\begin{Definition}\label{subdif} A vector $v\in\R^n$ is a subgradient of a convex function $f\colon\R^n\to(-\infty,\infty]$ at $\ox\in\dom(f)$ if it satisfies
the inequality
\begin{equation*}
f(x)\ge f(\ox)+\la v,x-\ox\ra \;\mbox{ for all }\;x\in\R^n.
\end{equation*}
The collection of subgradients is called the subdifferential of $f$ at $\ox$ and is denoted by $\partial f(\ox)$.
\end{Definition}

Directly from the definition  we have the {\em subdifferential Fermat rule}:
\begin{equation}\label{fermat}
\mbox{\rm $f$ has an absolute minimum at $\ox\in\dom(f)$ if and only if $0\in\partial f(\ox)$.}
\end{equation}

\begin{figure}[!ht]
\begin{center}
\includegraphics[width=9cm]{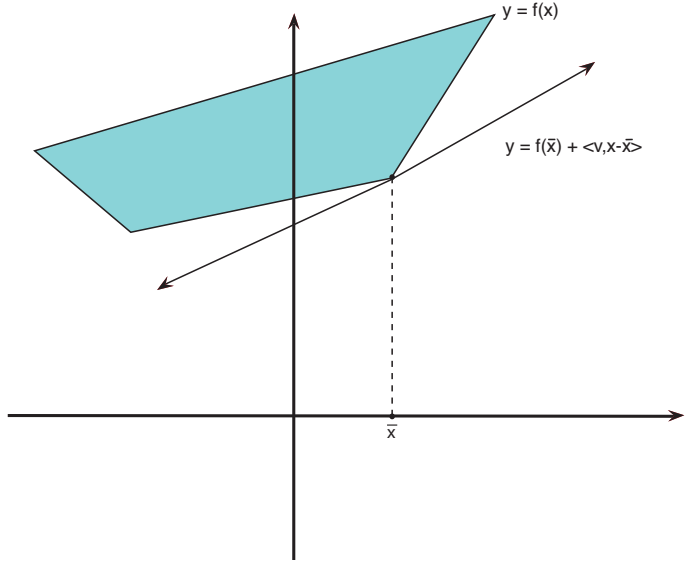}
\caption{Subgradient.}
\label{fig:1-sub}
\end{center}
\end{figure}

Recall that $f\colon\R^n\to(-\infty,\infty]$ is (Fr\'echet) {\em differentiable} at $\ox\in\dom(f)$ if there is a vector $v\in\R^n$ for which we have
\begin{equation*}
\lim_{x\to\ox}\frac{f(x)-f(\ox)-\la v,x-\ox\ra}{\|x-\ox\|}=0.
\end{equation*}
In this case the vector $v$ is unique, is known as the {\em gradient} of $f$ at $\ox$, and is denoted by $\nabla f(\ox)$. The next proposition shows that the subdifferential of a convex function at a given point reduces to its gradient at this point when the function is differentiable.

\begin{Proposition}\label{dd1} Let $f\colon\R^n\to(-\infty,\infty]$ is convex and differentiable at $\ox\in\dom(f)$. Then
\begin{equation}\label{d1}
\la\nabla f(\ox),x-\ox\ra\le f(x)-f(\ox)\;\mbox{ for all }\;x\in\R^n\;\mbox{ with }\;\partial f(\ox)=\{\nabla f(\ox)\}.
\end{equation}
\end{Proposition}
{\bf Proof.} Since $f$ is differentiable at $\ox$, for any $\epsilon>0$ there exists $\gamma>0$ such that
\begin{equation*}\label{d2}
-\epsilon\|x-\ox\|\le f(x)-f(\ox)-\la\nabla f(\ox),x-\ox\ra\le\epsilon\|x-\ox\|\;\mbox{\rm whenever }\;\|x-\ox\|\leq\gamma.
\end{equation*}
Define further the function
\begin{equation*}
\psi(x):=f(x)-f(\ox)-\la\nabla f(\ox),x-\ox\ra+\epsilon\|x-\ox\|
\end{equation*}
for which $\psi(x)\ge\psi(\ox)=0$ whenever $x\in\B(\ox;\gamma)$. It follows from the convexity of $\psi$ that $\psi(x)\ge\psi(\ox)$ when $x\in\R^n$, and thus
\begin{equation*}
\la\nabla f(\ox),x-\ox\ra\le f(x)-f(\ox)+\epsilon\|x-\ox\|\;\mbox{ for all }\;x\in\R^n.
\end{equation*}
Letting now $\epsilon\to 0^+$ gives us the inequality in (\ref{d1}) and shows that $\nabla f(\ox)\in\partial f(\ox)$.

To verify the remaining part of (\ref{d1}), pick any $v\in\partial f(\ox)$ and observe that
\begin{equation*}
\la v,x-\ox\ra\le f(x)-f(\ox)\;\mbox{\rm for all }\;x\in\R^n.
\end{equation*}
From the differentiability of $f$ at $\ox$ we have
\begin{equation*}
\la v-\nabla f(\ox),x-\ox\ra\le\epsilon\|x-\ox\|\;\mbox{ whenever }\;\|x-\ox\|\leq\gamma,
\end{equation*}
and so $\|v-\nabla f(\ox)\|\le\epsilon$. This yields $v=\nabla f(\ox)$ since $\epsilon>0$ is arbitrary and thus justifies the claimed relationship $\partial f(\ox)=\{\nabla f(\ox)\}$. $\h$

The following simple example calculates the subdifferential of the Euclidean norm function directly from the subdifferential definition.

\begin{Example} {\rm For the Euclidean norm function $p(x):=\|x\|$ we have
\begin{equation*}
\partial p(x)=\begin{cases}
\B &\text{if }\;x=0, \\
\Big\{\dfrac{x}{\|x\|}\Big\}& \text{otherwise},
\end{cases}
\end{equation*}
where $\B$ stands for the closed unit ball of $\R^n$. To verify this, we observe $\nabla p(x)=\dfrac{x}{\|x\|}$ for $x\ne 0$ due to the differentiability of $p(x)$ at nonzero points. Consider the case where $x=0$ and use the definition describe $v\in\partial p(0)$ as
\begin{equation*}
\la v,x\ra=\la v,x-0\ra\le p(x)-p(0)=\|x\|\;\mbox{\rm for all }\;x\in\R^n.
\end{equation*}
For $x=v$ therein we get $\la v,v\ra\le\|v\|$, which shows that $\|v\|\le 1$, i.e., $v\in\B$. Conversely, picking $v\in\B$ and employing the Cauchy-Schwarz inequality tell us that
\begin{equation*}
\la v,x-0\ra=\la v,x\ra\le\|v\|\cdot\|x\|\le\|x\|=p(x)-p(0)\;\mbox{\rm for all }\;x\in\R^n,
\end{equation*}
i.e., $v\in\partial p(0)$. Thus we arrive at the equality $\partial p(0)=\B$.}
\end{Example}

We conclude this section by the useful description of the relative interior of the graph of a convex function via the relative interior of its domain; see, e.g., \cite[Proposition~1.1.9]{HU1}.

\begin{Proposition}\label{re epi} For a convex function $f\colon\R^n\to(-\infty,\infty]$ we have the representation
\begin{equation*}
\ri(\epi(f))=\big\{(x,\lambda)\;\big|\;x\in\ri(\dom(f)),\;\lambda>f(x)\big\}.
\end{equation*}
\end{Proposition}

\section{Convex Separation}
\setcounter{equation}{0}

Separation theorems for convex sets, which go back to Minkowski \cite{min}, have been well recognized among the most fundamental geometric tools of convex analysis. In this section we formulate and give simple proofs of several separation results needed in what follows under the weakest assumptions in finite dimensions. Let us begin with \emph{strict separation} of a closed convex set and a point outside the set.

\begin{Proposition}\label{sep 2.1}
Let $\Omega$ be a nonempty closed convex set, and let $\bar{x}\notin\Omega$. Then there exists a nonzero vector $v\in\R^n$ such that
\begin{equation*}
\sup\{\langle v,x\rangle\;|\;x\in\Omega\}<\langle v,\bar{x}\rangle.
\end{equation*}
\end{Proposition}
{\bf Proof.} Denote $\bar{\omega}:=\pi(\bar{x};\Omega)$, $v:=\bar{x}-\bar{\omega}$, and fix $x\in\O$. Proposition~\ref{angle prop} gives us
$$
\langle v,x-\bar{\omega}\rangle=\langle\bar{x}-\bar{\omega},x-\bar{\omega}\rangle\le 0,
$$
which shows that
\begin{equation*}
\langle v,x-\bar\omega\rangle=\la v,x-\bar x+\bar x-\bar\omega\ra=\la v,x-\bar x+v\ra\le 0.
\end{equation*}
The last inequality therein yields
\begin{equation*}
\la v,x\ra\le\la v,\ox\ra-\|v\|^2,
\end{equation*}
which implies in turn that
$$
\sup\{\langle v,x\rangle\;|\;x\in\Omega\}< \langle v,\bar{x}\rangle
$$
and thus completes the proof of the proposition. $\h$

\begin{Remark}\label{rm1} {\rm It is easy to see that the closure $\Bar\O$ of a convex set $\O$ is convex. If $\O\subset\R^n$ is a nonempty convex set with $\ox\notin\Bar{\Omega}$, then applying Proposition~\ref{sep 2.1} to the convex set $\Bar\O$ gives us a nonzero vector $v\in\R^n$ such that
\begin{equation*}
\sup\{\langle v,x \rangle\;|\;x\in\Omega\}\le\sup\{\langle v,x\rangle\;|\;x\in\Bar\Omega\}<\la v,\ox\ra.
\end{equation*}}
\end{Remark}
\begin{figure}[!ht]
\begin{center}
\includegraphics[width=9cm]{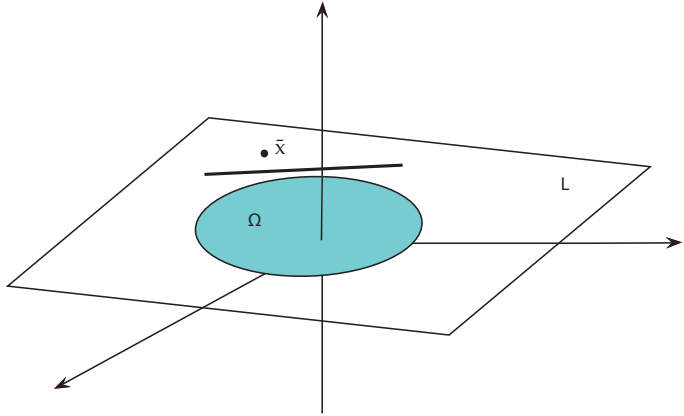}
 \caption{Separation in a subspace.}
\label{fig:1-sep}
\end{center}
\end{figure}

The next proposition justifies a strict separation property in a {\em subspace} of $\R^n$.

\begin{Proposition}\label{ss} Let $L$ be a subspace of $\R^n$, and let $\Omega\subset L$ be a nonempty convex set with $\bar{x}\in L$ and $\bar{x}\not\in\Bar{\Omega}$. Then there exists $v\in L$, $v\ne 0$, such that
$$
\sup\{\langle v,x\rangle\;|\;x\in\Omega\}<\langle v,\bar{x}\rangle.
$$
\end{Proposition}
{\bf Proof.} Employing Remark~\ref{rm1} gives us a vector $w\in\R^n$ such that
\begin{equation*}
\sup\{\langle w,x\rangle\;|\;x\in\Omega\}<\langle w,\bar{x}\rangle.
\end{equation*}
It is well known that $\R^n$ can be represented as the direct sum $\R^n=L\oplus L^{\perp}$, where
$$
L^{\perp}:=\{u\in\R^n\;|\;\langle u,x\rangle=0\;\text{ for all }\;x\in L\}.
$$
Thus $w=u+v$ with $u\in L^{\perp}$ and $v\in L$. This yields $\la u,x\ra=0$ for any $x\in\O\subset L$ and
\begin{eqnarray*}
\begin{array}{ll}
\langle v,x\rangle=\langle u,x\rangle+\langle v,x\rangle=\langle u+v,x\rangle=\langle w,x\rangle\le\sup\{\langle w,x\rangle\;|\;x\in\Omega\}\\
<\langle w,\bar{x}\rangle=\langle u+v,\bar{x}\rangle=\langle u,\bar{x}\rangle+\langle v,\bar{x}\rangle=\langle v,\bar{x}\rangle,
\end{array}
\end{eqnarray*}
which shows that $\sup\{\langle v,x \rangle\;|\;x\in\Omega\}<\langle v,\bar{x}\rangle$ with  $v\neq 0$. $\h$

\begin{figure}[!ht]
\begin{center}
\includegraphics[width=9cm]{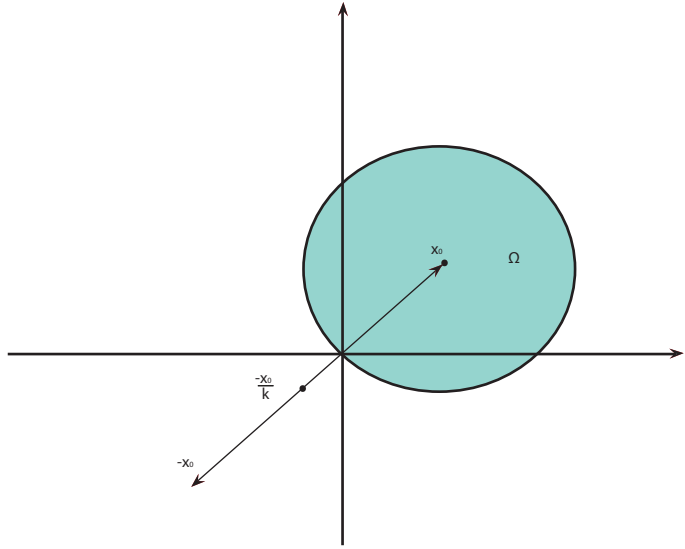}
\caption{Illustration of the proof of Lemma~\ref{lm3.1}.}
\label{fig:1-ill}
\end{center}
\end{figure}

\begin{Lemma}\label{lm3.1} Let $\O\subset\R^n$ be a nonempty convex set, and let $0\in\Bar{\Omega}\setminus\mbox{\rm ri}(\Omega)$. Then $\aff(\O)$ is a subspace of $\R^n$, and there is a sequence $\{x_k\}\subset\aff(\O)$ with $x_k\notin\Bar\O$ and $x_k\to 0$ as $k\to\infty$.
\end{Lemma}
{\bf Proof.} Since $\ri(\Omega)\ne\emp$ by Theorem~\ref{ri 1.72}(i) and $0\in\Bar{\Omega}\setminus\mbox{\rm ri}(\Omega)$, we find $x_0\in\ri(\Omega)$ and conclude that $-tx_0\notin\Bar{\Omega}$ for all $t>0$. Indeed, suppose by contradiction that $-tx_0\in\Bar{\Omega}$ for some $t>0$ and then deduce from Theorem~\ref{ri 1.72}(ii) that
\begin{equation*}
0=\frac{t}{1+t}x_0+\frac{1}{1+t}(-tx_0)\in\mbox{\rm ri}(\Omega),
\end{equation*}
which contradicts $0\notin\ri(\Omega)$. Letting now $x_k:=-\frac{x_0}{k}$ implies that $x_k\notin\Bar{\Omega}$ for every $k$ and $x_k\to 0$ as $k\to\infty$. Furthermore, we have
$$
0\in\Bar{\O}\subset\Bar{\aff(\O)}=\aff(\O)
$$
by the closedness of $\aff(\O)$ due to Proposition~\ref{closed property}(ii). This shows that $\aff(\O)$ is a subspace and that $x_k\in\aff(\O)$ for all $k\in\N$. $\h$

We continue with another important separation property known as \emph{proper separation}.

\begin{Definition} It is said that two nonempty convex sets $\O_1$ and $\O_2$ are properly separated if there exists a nonzero vector $v\in\R^n$ such that
\begin{equation*}
\sup\{\la v, x\ra\;|\; x\in \O_1\}\le\inf\{\la v,y\ra\;|\;y\in\O_2\},\quad\inf\{\la v, x\ra\;|\;x\in\O_1\}<\sup\{\la v,y\ra\;|\;y\in\O_2\}.
\end{equation*}
\end{Definition}

\begin{Lemma}\label{prop3.6} Let $\Omega$ be a nonempty convex set in $\R^n$. Then $0\notin\ri(\Omega)$ if and only if the sets $\Omega$ and $\{0\}$ are properly separated, i.e., there is $v\in\R^n$, $v\ne 0$, such that
$$
\sup\{\langle v,x\rangle\;|\; x\in\Omega\}\le 0,\quad\inf\{\langle v,x\rangle\;|\;x\in \Omega\}<0.
$$
\end{Lemma}
{\bf Proof.}  We split the proof into the following two cases.\\[1ex]
\emph{Case 1:} $0\not\in\Bar{\Omega}$. It follows from Remark~\ref{rm1} with $\ox=0$ that there exists $v\ne 0$ such that
$$
\sup\{\langle v,x\rangle\;|\;x\in\Omega\}<\la v,\ox\ra=0,
$$
and thus the sets $\Omega$ and $\{0\}$ are properly separated.\\[1ex]
\emph{Case 2:} $0\in\Bar{\Omega}\setminus\text{\rm ri}(\Omega)$. Letting $L:=\mbox{\rm aff}(\O)$ and employing Lemma~\ref{lm3.1} tell us that $L$ is a subspace of $\R^n$ and there is a sequence $\{x_k\}\subset L$ with $x_k\notin\Bar{\Omega}$ and $x_k\to 0$ as $k\to\infty$. By Proposition~\ref{ss} there is a sequence $\{v_k\}\subset L$ with $v_k\ne 0$ and
\begin{equation*}
\sup\{\langle v_k,x\rangle\;|\;x\in\Omega\}<\la v_k,x_k\ra,\quad k\in\N.
\end{equation*}
Denoting $w_k:=\frac{v_k}{\Vert v_k\Vert}$ shows that $\|w_k\|=1$ for all $k\in\N$ and
\begin{equation}\label{wk}
\langle w_k,x\rangle<\la w_k,x_k\ra\;\mbox{\rm for all }\;x\in\O.
\end{equation}
Letting $k\to\infty$ in (\ref{wk}) and supposing without loss of generality that $w_k \to v\in L$ with some $\|v\|=1$ along the whole sequence of $\{k\}$, we arrive at
\begin{equation*}
\sup\{\langle v,x\rangle\;|\;x\in\Omega\}\le 0
\end{equation*}
by taking into account that $|\la w_k,x_k\ra|\le\|w_k\|\cdot\|x_k\|=\|x_k\|\to 0$. To verify further
$$
\inf\{\langle v,x\rangle\;|\;x\in\Omega\}<0,
$$
it suffices to show that there is $x\in\Omega$ with $\langle v,x\rangle <0$. Suppose by contradiction that $\langle v,x\rangle\ge 0$ for all $x\in\Omega$ and deduce from $\sup\{\langle v,x\rangle\;|\;x\in\Omega\}\le 0$ that $\langle v,x\rangle=0$ for all $x\in\Omega$. Since $v\in L=\aff(\Omega)$, we get the representation
$$
v=\sum_{i=1}^{m}\lambda_i\omega_i\;\mbox{ with }\;\sum_{i=1}^{m}\lambda_i=1\;\mbox{ and }\;\omega_i\in\Omega\;\mbox{ for }\;i=1,\ldots,m,
$$
which readily implies the equalities
\begin{align*}
\Vert v\Vert^2=\langle v,v\rangle=\sum_{i=1}^{m}\lambda_i\langle v,\omega_i\rangle=0
\end{align*}
and so contradicts the condition $\|v\|=1$. This justifies the proper separation of $\O$ and $\{0\}$.

To verify the reverse statement of the lemma, assume that $\Omega$ and $\{0\}$ are properly separated and thus find $0\ne v\in\R^n$ such that
\begin{equation*}
\sup\{\langle v,x\rangle\;|\;x\in\Omega\}\le 0\;\mbox{ while }\;\langle v,\bar{x}\rangle<0\;\mbox{ for some }\;\bar{x}\in\Omega.
\end{equation*}
Suppose by contradiction that $0\in\ri(\Omega)$ and deduce from Proposition~\ref{prop 1.5} that
$$
0+t(0-\bar{x})=-t\bar{x}\in\Omega\;\mbox{\rm for some }\;t>0.
$$
This immediately implies the inequalities
$$
\langle v,-t\bar{x}\rangle\le\sup\{\langle v,x\rangle\;|\;x \in\Omega\}\le 0
$$
showing that $\langle v,\bar{x}\rangle\ge 0$. It is a contradiction, which verifies $0\notin\ri(\Omega)$. $\h$

Now we are ready to prove the {\em main separation theorem} in convex analysis.

\begin{Theorem}\label{propsep}
Let $\Omega_1$ and $\Omega_2$ be two nonempty convex subsets of $\R^n$. Then $\Omega_1$ and $\Omega_2$ are properly separated if and only if $\ri(\Omega_1)\cap \ri(\Omega_2)=\emp$.
\end{Theorem}
{\bf Proof.} Define $\Omega:=\Omega_1-\Omega_2$ and verify that $\ri(\Omega_1)\cap\ri(\Omega_2)=\emp$ if and only if
$$
0\notin\ri(\Omega_1-\Omega_2)=\ri(\Omega_1)-\ri(\Omega_2).
$$
To proceed, suppose first that $\ri(\Omega_1)\cap\ri(\Omega_2)=\emp$ and so get by Corollary~\ref{ri of set difference} that $0\not\in\ri(\Omega_1-\Omega_2)=\ri(\Omega)$. Then Lemma~\ref{prop3.6} tells us that the sets $\Omega$ and $\{0\}$ are properly separated. Thus there exist $v\in\R^n$ with $\langle v,x\rangle\le 0$ for all $x\in\Omega$ and also $y\in\Omega$ such that $\langle v,y\rangle<0$. For any $\omega_1\in\O_1$ and $\omega_2\in\O_2$ we have $x:=\omega_1-\omega_2\in\O$, and hence
\begin{align*}
\la v,\omega_1-\omega_2\ra=\la v,x\ra\le 0,
\end{align*}
which yields $\la v,\omega_1\ra\le\la v,\omega_2\ra$. Choose $\bar\omega_1\in\O_1$ and $\bar\omega_2\in\O_2$ such that $y=\bar\omega_1-\bar\omega_2$. Then
\begin{equation*}
\la v,\bar\omega_1-\bar\omega_2\ra=\la v,y\ra<0
\end{equation*}
telling us that $\la v,\bar\omega_1\ra<\la v,\bar\omega_2\ra$. Hence $\Omega_1$ and $\Omega_2$ are properly separated.

To prove next the converse implication, suppose that $\Omega_1$ and $\Omega_2$ are properly separated, which implies that the sets $\Omega=\Omega_1-\Omega_2$ and $\{0\}$ are properly separated as well. Employing Lemma~\ref{prop3.6} again provides the relationships
$$
0\notin\ri(\Omega)=\ri(\Omega_1-\Omega_2)=\ri(\Omega_1)-\ri(\Omega_2)\;\mbox{ and so }\;\ri(\Omega_1)\cap\ri(\Omega_2)=\emp,
$$
which thus complete the proof of the theorem. $\h$

Among various consequences of Theorem~\ref{propsep}, including those presented below, note the following relationships between the closure and relative interior operations on convex sets that seemingly have nothing to do with separation.

\begin{Corollary}\label{rc} {\bf (i)} If $\O\subset\R^n$ is convex, then $\ri(\Bar\O)=\ri(\O)$ and $\Bar{\ri(\O)}=\Bar\O$.\\
{\bf (ii)} If both sets $\O_1,\O_2\subset\R^n$ are convex, then we have $\ri(\O_1)=\ri(\O_2)$ provided that $\Bar{\O}_1=\Bar{\O}_2$.
\end{Corollary}
{\bf Proof.} {\bf(i)} Both equalities are trivial if $\O=\emp$. To verify the first equality in (i) when $\O\ne\emp$, observe that for any $x\in\R^n$ we have the equivalence
$$
\{x\}\;\mbox{\rm and }\O\;\mbox{\rm are properly separated }\Longleftrightarrow\{x\}\;\mbox{\rm and }\Bar{\O}\;\mbox{\rm are properly separated}.
$$
Indeed, the implication ``$\Longrightarrow$'' is obvious because  for any $v\in \R^n$
\begin{equation*}
\big[\sup\{\langle v,w \rangle\;|\;w\in\Omega\}\leq \la v,x\ra\big]\Longrightarrow \Big[\sup\{\langle v,w\rangle\;|\;w\in\Bar\Omega\}\leq \la v, x\ra\big],
\end{equation*}
which can be proved by a limiting argument. For the converse, suppose that $\{x\}$ and $\Bar{\O}$ are properly separated. Then there exists $v\in \R^n$, $v\neq 0$, such that
\begin{equation*}
\sup\{\la v, w\ra\;|\; w\in \O\}\le \la v, x\ra,\quad\inf\{\la v, w\ra\;|\;w\in\O\}<\la v,x\ra.
\end{equation*}
It follows that
\begin{equation*}
\sup\{\langle v,w \rangle\;|\;w\in\Omega\}\leq \sup\{\langle v,w\rangle\;|\;w\in\Bar\Omega\}\leq \la v, x\ra.
\end{equation*}
It remains to show that there exists $\bar{w}\in \O$ such that $\la v, \bar{w}\ra < \la v, x\ra$. If this is not the case, then $\la v, w\ra\geq \la v,x\ra$ for all $w\in \O$, which implies $\la v, w\ra=\la v,x\ra$ for all $w\in \O$. A simple limiting argument yields $\la v, w\ra=\la v,x\ra$ for all $w\in \Bar\O$. This is a contradiction because $\inf\{\langle v,w \rangle\;|\;w\in\Bar\Omega\}<\la v,x\ra$.

Defining further $\Theta:=\{x\}$, we get $\ri(\Theta)=\{x\}$ and deduce from Theorem~\ref{propsep} that
\begin{align*}
x\notin\ri(\Bar\O)&\Longleftrightarrow\ri(\Theta)\cap\ri(\Bar{\O})=\emp\\
&\Longleftrightarrow\{x\}\;\mbox{\rm and }\;\Bar{\O}\;\mbox{\rm are properly separated}\\
&\Longleftrightarrow\{x\}\;\mbox{\rm and }\O\;\mbox{\rm are properly separated}\\
&\Longleftrightarrow\ri(\Theta)\cap\ri(\O)=\emp\Longleftrightarrow x\notin\ri(\O).
\end{align*}
The second equality in (i) is a direct consequence of Theorem~\ref{ri 1.72}(ii).\\[1ex]
{\bf(ii)} If $\Bar{\O}_1=\Bar{\O}_2$, then $\ri(\Bar{\O}_1)=\ri(\Bar{\O}_1)$ and hence $\ri(\O_1)=\ri(\O_2)$ by (i). $\h$

\section{Normal Cone Intersection Rule}
\setcounter{equation}{0}

In this section we derive the central result of the geometric approach to convex subdifferential calculus, which provides a general {\em intersection rule} for the normal cone to convex sets. All the subsequent subdifferential results are consequences of this intersection rule.

Recall first the definition of the normal cone to a convex set.

\begin{Definition}\label{nc-def} Let $\O$ be a nonempty convex subset of $\R^n$. Then the normal cone to the set $\O$ at $\ox\in\O$ is defined by
\begin{equation*}
N(\ox;\O):=\{v\in\R^n\;|\;\la v,x-\ox\ra\le 0\;\mbox{ for all }\;x\in\O\}.
\end{equation*}
In the case where $\ox\notin\O$ we define $N(\ox;\O):=\emp$.
\end{Definition}

\begin{figure}[!ht]
\begin{center}
\includegraphics[width=10cm]{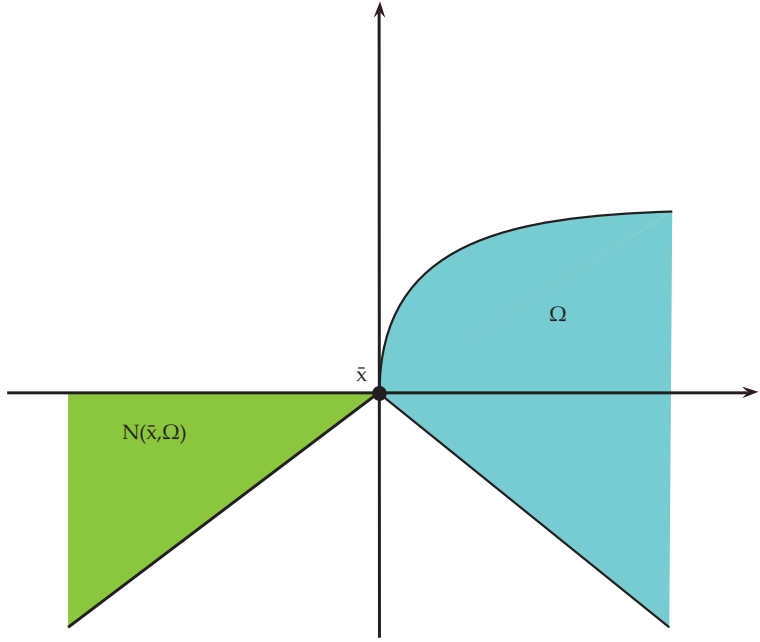}
 \caption{Normal cone.}
\label{fig:1-nc}
\end{center}
\end{figure}

It immediately follows from Definition~\ref{nc} that $N(\ox;\O)$ is a closed and convex cone, which reduces to $\{0\}$ if $\ox\in{\rm int}(\O)$. A remarkable property of the normal cone to a convex  set $\O$ in finite dimensions is that $N(\ox;\O)\ne\{0\}$ if and only if $\ox$ is boundary point of $\O$; see, e.g., \cite[Corollary~2.14]{bmn}. This is the classical {\em supporting hyperplane theorem}, which can be easily derived  by the limiting procedure from Theorem~\ref{propsep}.

Before deriving our major intersection result on the representation of the normal cone to finitely many convex sets, let us present a useful lemma on the relative interior of set intersections, which is also based on convex separation.

\begin{Lemma}\label{ri} Let $\O_i\subset\R^n$ for $i=1,\ldots,m$ with $m\ge 2$ be convex subsets of $\R^n$ such that
\begin{eqnarray}\label{ri-qc}
\bigcap_{i=1}^m\ri(\O_i)\ne\emp.
\end{eqnarray}
Then we have the representation
\begin{eqnarray}\label{ri-inter}
\ri\Big(\bigcap_{i=1}^m\O_i\Big)=\bigcap_{i=1}^m\ri(\O_i).
\end{eqnarray}
\end{Lemma}
{\bf Proof.} We first verify this result for $m=2$. Pick $x\in\ri(\O_1)\cap\ri(\O_2)$ and find $\gamma>0$ with
\begin{equation*}
\B(x;\gamma)\cap\aff(\O_1)\subset\O_1\;\mbox{\rm and }\;\B(x;\gamma)\cap\aff(\O_2)\subset\O_2,
\end{equation*}
which implies therefore that
\begin{equation*}
\B(x;\gamma)\cap[\aff(\O_1)\cap\aff(\O_2)]\subset\O_1\cap\O_2.
\end{equation*}
It is easy to see that $\aff(\O_1\cap\O_2)\subset\aff(\O_1)\cap\aff(\O_2)$, and hence
\begin{equation*}
\B(x;\gamma)\cap\aff(\O_1\cap\O_2)\subset\O_1\cap\O_2.
\end{equation*}
Thus we get $x\in\ri(\O_1\cap\O_2)$, which justifies that $\ri(\O_1)\cap\ri(\O_2)\subset\ri(\O_1\cap\O_2)$.

To verify the opposite inclusion in (\ref{ri-inter}) for $m=2$, observe that
\begin{eqnarray}\label{cl}
\Bar{\O_1\cap\O_2}=\Bar{\O}_1\cap\Bar{\O}_2
\end{eqnarray}
for any convex sets $\O_1,\O_2$ with $\ri(\O_1)\cap\ri(\O_2)\ne\emp$. Indeed, pick $x\in\Bar{\O}_1\cap\Bar{\O}_2$, $\ox\in\ri(\O_1)\cap\ri(\O_2)$ and observe that $x_k:=k^{-1}\ox+(1-k^{-1})x\to x$ as $k\to\infty$. Then Theorem~\ref{ri 1.72}(ii) tells us that $x_k\in\O_1\cap\O_2$ for large $k\in\N$ and hence $x\in \Bar{\O_1\cap\O_2}$, which justifies the inclusion ``$\supset$" in (\ref{cl}). The inclusion ``$\subset$" therein obviously holds even for nonconvex sets. Now using (\ref{cl}) and the second equality in Corollary~\ref{rc}(i) gives us
\begin{equation*}
\Bar{\ri(\O_1)\cap\ri(\O_2)}=\Bar{\ri(\O_1)}\cap\Bar{\ri(\O_2)}=\Bar{\O_1}\cap\Bar{\O_2}=\Bar{\O_1\cap\O_2}.
\end{equation*}
Then we have the equality
\begin{equation*}
\ri(\Bar{\ri(\O_1)\cap\ri(\O_2)})=\ri(\Bar{\O_1\cap\O_2}),
\end{equation*}
and thus conclude by Corollary~\ref{rc}(ii) that
\begin{equation*}
\ri(\ri(\O_1)\cap\ri(\O_2))=\ri(\O_1\cap\O_2)\subset\ri(\O_1)\cap\ri(\O_2),
\end{equation*}
which justify representation (\ref{ri-inter}) for $m=2$.

To verify (\ref{ri-inter}) under the validity of (\ref{ri-qc}) in the general case of $m>2$, we proceed by induction with taking into account that the result has been established for two sets and assuming that it holds for $m-1$ sets. Considering $m$ sets $\O_i$, represent their intersection as
\begin{eqnarray}\label{inter}
\disp\bigcap_{i=1}^m\O_i=\O\cap\O_m\;\mbox{ with }\;\O:=\bigcap_{i=1}^{m-1}\O_i.
\end{eqnarray}
Then we have $\ri(\O)\cap\ri(\O_m)=\cap_{i=1}^m\ri(\O_i)\ne\emp$ by the imposed condition in (\ref{ri-qc}) and the induction assumption on the validity of  (\ref{ri-inter}) for $m-1$ sets. This allows us to employ the obtained result for the two sets $\O$ and $\O_m$ and thus arrive at the desired conclusion  (\ref{ri-inter}) for the $m$ sets $\O_1,\ldots,\O_m$ under consideration. $\h$

Now we are ready to derive the underlying formula for the representation of the normal cone to intersections of finitely many convex sets. Note that the proof of this result and of the subsequent calculus rules for functions and set-valued mappings mainly follow the geometric pattern of variational analysis as in \cite{m-book1}. The specific features of convexity and the usage of convex separation instead of the extremal principle allow us to essentially simplify the proof and to {\em avoid the closedness} requirement on sets and the corresponding {\em lower semicontinuity} assumptions on functions in subdifferential calculus rules. Furthermore, we show below that the developed geometric approach works in the convex setting under the {\em relative interior qualification conditions}, which are well-recognized in finite-dimensional convex analysis and occur to be {\em weaker} than the basic/normal qualifications employed in \cite{m-book1,bmn}; see. e.g., Corollary~\ref{bqc1} and Example~\ref{ex-qc} below.

\begin{figure}[!ht]
\begin{center}
\includegraphics[width=8cm]{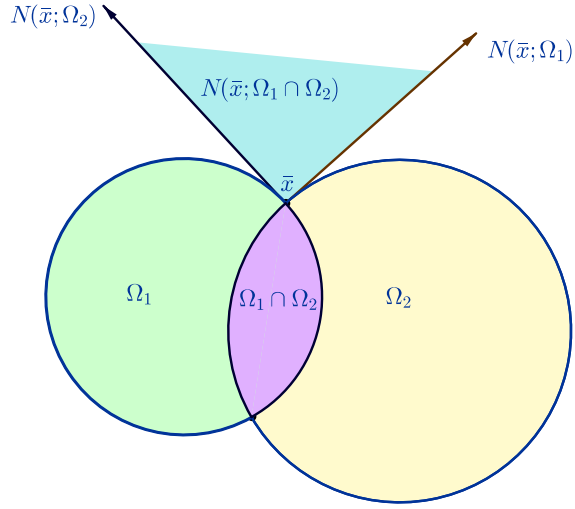}
\caption{Intersection rule.}
\label{fig:1-int1}
\end{center}
\end{figure}

\begin{Theorem}\label{nc} Let $\O_1,\ldots,\O_m\subset\R^n$ be convex sets satisfying the relative interior condition
\begin{eqnarray}\label{qc1}
\bigcap_{i=1}^m\ri(\O_i)\ne\emp,
\end{eqnarray}
where $m\ge 2$. Then we have the intersection rule
\begin{equation}\label{sum2}
N\Big(\bar{x};\bigcap_{i=1}^m\O_i)\Big)=\disp\sum_{i=1}^m N(\bar{x};\Omega_i)\;\mbox{\rm for all }\bar{x}\in\bigcap_{i=1}^m\O_i.
\end{equation}
\end{Theorem}
{\bf Proof.} Proceeding by induction, let us first prove the statement of the theorem for the case of $m=2$. Since the inclusion ``$\supset$" in (\ref{sum2}) trivially holds even without imposing \eqref{qc1}, the real task is to verify the opposite inclusion therein. Fixing $\ox\in\Omega_1\cap\Omega_2$ and $v\in N(\bar x;\Omega_1\cap\Omega_2)$, we get by the normal cone definition that
\begin{equation*}
\langle v,x-\bar x\rangle\le 0\;\mbox{ for all }\;x\in\Omega_1\cap\Omega_2.
\end{equation*}
Denote $\Theta_1:=\Omega_1\times [0,\infty)$ and $\Theta_2:=\{(x,\lambda)\;|\;x\in\Omega_2,\;\lambda\le\langle v,x-\bar x\rangle\}$. It follows from Proposition~\ref{re epi} that $\ri(\Theta_1)=\ri(\O_1)\times(0,\infty)$ and
\begin{equation*}
\ri(\Theta_2)=\big\{(x,\lambda)\;\big|\;x\in\ri(\O_2),\;\lambda<\la v,x-\ox\ra\big\}.
\end{equation*}
Arguing by contradiction, it is easy to check that $\ri(\Theta_1)\cap\ri(\Theta_2)=\emp$. Then applying Theorem~\ref{propsep} to these convex sets in $\R^{n+1}$ gives us $0\ne(w,\gamma)\in\R^n\times\R$ such that
\begin{equation}\label{separation}
\la w,x\ra+\lambda_1\gamma\le\la w,y\ra+\lambda_2\gamma\;\mbox{ for all }\;(x,\lambda_1)\in \Theta_1,\;(y,\lambda_2)\in \Theta_2.
\end{equation}
Moreover, there are $(\tilde{x},\tilde{\lambda}_1)\in\Theta_1$ and $(\tilde{y},\tilde{\lambda}_2)\in\Theta_2$ satisfying
\begin{equation*}
\la w,\tilde{x}\ra+\tilde{\lambda}_1\gamma<\la w,\tilde{y}\ra+\tilde{\lambda}_2\gamma.
\end{equation*}
Observe that $\gamma\le 0$ since otherwise we can get a contradiction by employing \eqref{separation} with $(\ox,k)\in\Theta_1$ as $k>0$ and $(\ox,0)\in\Theta_2$. Let us now show by using (\ref{qc1}) that $\gamma<0$. Again arguing by contradiction, suppose that $\gamma=0$ and then get
\begin{align*}
\la w,x\ra\le\la w,y\ra\;\mbox{ for all }\; x\in\O_1,\;y\in\O_2\;\mbox{ and }\;\la w,\tilde{x}\ra<\la w,\tilde{y}\ra\;\mbox{ with }\;\tilde{x}\in\O_1,\;\tilde{y}\in\O_2.
\end{align*}
This means the proper separation of the sets $\O_1$ and $\O_2$, which tells us by Theorem~\ref{propsep} that $\ri(\O_1)\cap\ri(\O_2)=\emp$. The obtained contradiction verifies the claim of $\gamma<0$.

To proceed further, denote $\mu:=-\gamma>0$ and deduce from (\ref{separation}), by taking into account that
$(x,0)\in\Theta_1$ when $x\in\Omega_1$ and and that $(\ox,0)\in\Theta_2$, the inequality
\begin{equation*}
\la w,x\ra\le\la w,\ox\ra\;\mbox{ for all }\;x\in\Omega_1.
\end{equation*}
This yields $w\in N(\ox;\Omega_1)$ and hence $\dfrac{w}{\mu}\in N(\ox;\Omega_1)$. Moreover, we get from \eqref{separation}, due to $(\ox,0)\in \Theta_1$ and $(y,\al)\in\Theta_2$ for all $y\in\O_2$ with $\al=\la v,y-\ox\ra$, that
\begin{equation*}
\big\la w,\ox\big\ra\le\big\la w,y\big\ra+\gamma\la v,y-\ox\ra\;\mbox{ whenever }\;y\in\Omega_2.
\end{equation*}
Dividing both sides therein by $\gamma$, we arrive at the relationship
\begin{equation*}
\Big\la\dfrac{w}{\gamma}+v,y-\ox\Big\ra\le 0\;\mbox{ for all }\;y\in\Omega_2,
\end{equation*}
and thus $\dfrac{w}{\gamma}+v=-\dfrac{w}{\mu}+v\in N(\ox;\Omega_2)$. This gives us
\begin{equation*}
v\in\dfrac{w}{\mu}+N(\ox;\Omega_2)\subset N(\ox;\O_1)+N(\ox;\O_2)
\end{equation*}
completing therefore the proof of (\ref{sum2}) in the case of $m=2$.

Considering now the case of intersections for any finite number of sets, suppose by induction that the intersection rule (\ref{sum2}) holds under (\ref{qc1}) for $m-1$ sets and verify that it continues to hold for the intersection of $m>2$ sets $\bigcap_{i=1}^m\O_i$. Represent the latter intersection as $\O\cap\O_m$ with $\O:=\bigcap_{i=1}^{m-1}\O_i$, we get from the imposed relative interior condition (\ref{qc1}) and Lemma~\ref{ri} that
\begin{eqnarray*}
\ri(\O)\cap\ri(\O_m)=\bigcap_{i=1}^m\ri(\O_i)\neq\emp.
\end{eqnarray*}
Applying the intersection rule (\ref{sum2}) to the two sets $\O\cap\O_m$ and then employing the induction assumption for $m-1$ sets give us the equalities
\begin{eqnarray*}
N\Big(\ox;\bigcap_{i=1}^m\O_i\big)=N(\ox;\O\cap\O_m)=N(\ox;\O)+N(\ox;\O_m)=\disp\sum_{i=1}^m N(\ox;\O_i),
\end{eqnarray*}
which thus justify (\ref{sum2}) for $m$ sets and thus completes the proof of the theorem. $\h$

It is not difficult to observe the relative interior assumption (\ref{qc1}) is {\em essential} for the validity of the intersection rule (\ref{sum2}) as illustrated by the following example.

\begin{Example}{\rm Define the two convex sets on the plane by
\begin{equation*}
\O_1:=\{(x,\lambda)\in\R^2\;|\;\lambda\ge x^2\}\;\mbox{\rm and }\;\O_2:=\{(x,\lambda)\in \R^2\;|\;\lambda\le-x^2\}.
\end{equation*}
Then for $\ox=(0,0)\in\O_1\cap\O_2$ we have
\begin{equation*}
N(\ox;\O_1)=\{0\}\times(-\infty,0],\quad N(\ox;\O_2)=\{0\}\times[0,\infty),\;\mbox{\rm and }\;N(\ox;\O_1\cap\O_2)=\R^2.
\end{equation*}
Thus $N(\ox;\O_1)+N(\ox;\O_2)=\{0\}\times\R\ne N(\ox;\O_1\cap\O_2)$, i.e., the intersection rule (\ref{sum2}) fails. It does not contradict Theorem~\ref{nc}, since $\mbox{\rm ri}(\O_1)\cap \mbox{\rm ri}(\O_2)=\emp$ and so the relative interior qualification condition (\ref{qc1}) does not hold in this case.}
\end{Example}

\begin{figure}[!ht]
\begin{center}
\includegraphics[width=7cm]{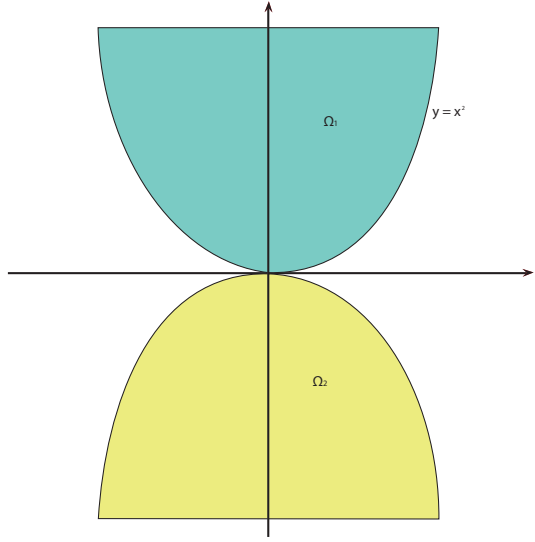}
\caption{Illustration of the relative interior condition.}
\label{fig:1a}
\end{center}
\end{figure}

Finally, we compare the intersection rule of Theorem~\ref{nc} derived under the relative interior qualification condition (\ref{qc1}) with the corresponding result obtained in \cite[Corollary~2.16]{bmn} for $m=2$ under the so-called {\em basic/normal qualification condition}
\begin{equation}\label{bqc}
N(\ox;\O_1)\cap[-N(\ox;\O_2)]=\{0\}
\end{equation}
introduced and applied earlier for deriving the intersection rule and related calculus results in nonconvex variational analysis; see, e.g., \cite{m-book1,rw} and the references therein. Let us first show that (\ref{bqc}) yields (\ref{qc1}) in the general convex setting.

\begin{Corollary}\label{bqc1} Let $\O_1,\O_2\subset\R^n$ be convex sets satisfying the basic qualification condition {\rm(\ref{bqc})} at some $\ox\in\O_1\cap\O_2$. Then we have
\begin{equation}\label{qc2}
\mbox{\rm ri}(\O_1)\cap\mbox{\rm ri}(\O_2)\ne\emp,
\end{equation}
and so the intersection rule {\rm(\ref{sum2})} holds for these sets at any $\ox\in \O_1\cap\O_2$.
\end{Corollary}
{\bf Proof}. Arguing by contradiction, suppose that $\mbox{\rm ri}(\O_1)\cap\mbox{\rm ri}(\O_2)=\emp$. Then the sets $\O_1,\O_2$ are properly separated by
Theorem~\ref{propsep}, and so there is $v\ne 0$ such that
\begin{equation*}
\la v,x\ra\le\la v,y\ra\;\mbox{\rm for all }\;x\in\O_1,\;y\in\O_2.
\end{equation*}
Since $\ox\in\O_2$, we have $\la v,x-\ox\ra\le 0$ for all $x\in\O_1$. Hence $v\in N(\ox;\O_1)$ and similarly $-v\in N(\ox;\O_2)$. Thus $0\ne v\in N(\ox;\O_1)\cap[-N(\ox;\O_2)]$, which contradicts (\ref{bqc}). $\h$

The next example demonstrates that (\ref{qc2}) may be {\em strictly weaker} then (\ref{bqc}).

\begin{Example}\label{ex-qc}{\rm Consider the two convex sets on the plane defined by $\O_1:=\R\times\{0\}$ and $\O_2:=(-\infty,0]\times\{0\}$. We obviously
get that condition (\ref{qc2}) is satisfied ensuring thus the validity of the intersection rule by Theorem~\ref{nc}. On the other hand, it follows for
that
\begin{equation*}
N(\ox;\O_1)=\{0\}\times\R\;\mbox{\rm and }\;N(\ox;\O_2)=[0,\infty)\times\R\;\mbox{ with }\;\ox=(0,0),
\end{equation*}
i.e., the other qualification condition (\ref{bqc}) fails, which shows that the result of \cite[Corollary~2.16]{bmn} is not applicable in this case.}
\end{Example}

\section{Subdifferential Sum Rule and Existence of Subgradients}
\setcounter{equation}{0}

The main goal of this section is to derive from the geometric intersection rule of Theorem~\ref{nc} the subdifferential sum rule for convex extended-real-valued functions under the least restrictive relative interior qualification condition. Then we deduce from it a mild condition ensuring the existence of subgradients for general convex functions.

Prior to this, let us recall well-known relationships between normals to convex sets and subgradients of convex functions used in what follows.

\begin{Proposition}\label{snr} {\bf (i)} Let $\O\subset\R^n$  be a nonempty convex set, and let $\delta(x;\O)=\dd_\O(x)$ be its indicator function equal to $0$ when $x\in\O$ and to $\infty$ otherwise. Then we have
$$
\partial\delta(\ox;\O)=N(\ox;\O)\;\mbox{ for any }\;\ox\in\O.
$$

{\bf(ii)} Let $f:\R^n\to(-\infty,\infty]$ be a convex function, and let $\ox\in\mbox{\rm dom}(f)$. Then we have
\begin{equation*}
\partial f(\ox)=\{v\in\R^n\;|\;(v,-1)\in N((\ox,f(\ox));\epi(f))\}.
\end{equation*}
\end{Proposition}
{\bf Proof.} {\bf (i)} It follows directly from the definitions of the subdifferential, normal cone, and the set indicator function.\\[1ex]
{\bf (ii)} Fix any subgradient $v\in\partial f(\ox)$ and then get from Definition~\ref{subdif} that
\begin{equation}\label{si}
\la v,x-\ox\ra\le f(x)-f(\ox)\;\mbox{\rm for all }\;x\in\R^n.
\end{equation}
To show that $(v,-1)\in N((\ox,f(\ox));\epi(f))$, fix any $(x,\lambda)\in\epi(f)$ and observe that due to $\lambda\ge f(x)$ we have the relationships
\begin{align*}
&\la(v,-1),(x,\lambda)-(\ox,f(\ox))\ra=\la v, x-\ox\ra+(-1)(\lambda-f(\ox))\\
&=\la v,x-\ox\ra-(\lambda-f(\ox))\le\la v,x-\ox\ra-(f(x)-f(\ox))\le 0,
\end{align*}
where the the last inequality holds by \eqref{si}. To verify the opposite inclusion in (ii), take $(v,-1)\in N((\ox,f(\ox));\epi(f))$ and fix any $x\in\mbox{\rm dom}(f)$. Then $(x,f(x))\in\epi(f)$ and hence
\begin{equation*}
\la (v,-1),(x, f(x))-(\ox, f(\ox))\ra\le 0,
\end{equation*}
which in turn implies the inequality
\begin{equation*}
\la v,x-\ox\ra-(f(x)-f(\ox))\le 0.
\end{equation*}
Thus $v\in\partial f(\ox)$, which completes the proof of the proposition. $\h$

Now we are ready to deduce the following subdifferential sum rule for function from the intersection rule of Theorem~\ref{nc} for normals to sets.

\begin{Theorem}\label{sr}
Let $f_i\colon\R^n\to(-\infty,\infty]$, $i=1,\ldots,m$, be extended-real-valued convex functions satisfying the relative interior qualification condition
\begin{equation}\label{riq}
\bigcap_{i=1}^m\ri\Big(\dom(f_i)\Big)\ne\emp,
\end{equation}
where $m\ge 2$. Then for all $\ox\in\bigcap_{i=1}^m\dom(f_i)$ we have the sum rule
\begin{equation}\label{ssr}
\partial\Big(\sum_{i=1}^m f_i\Big)(\ox)=\sum_{i=1}^m\partial f_i(\ox).
\end{equation}
\end{Theorem}
{\bf Proof.} Observing that the inclusion ``$\supset$" in (\ref{ssr}) directly follows from the subdifferential definition, we proceed with the proof of the opposite inclusion. Consider first the case of $m=2$ and pick any $v\in\partial(f_1+f_2)(\bar{x})$. Then we have
\begin{equation}\label{by def sub}
\langle v,x-\bar{x}\rangle\le(f_1+f_2)(x)-(f_1+f_2)(\bar{x})\;\mbox{\rm for all }\;x\in\R^n.
\end{equation}
Define the following convex subsets of $\R^{n+2}$ by
\begin{align*}
&\Omega_1:=\{(x,\lambda_1,\lambda_2)\in\R^n\times\R\times\R\;|\;\lambda_1\ge f_1(x)\}=\epi(f_1)\times\R,\\
&\Omega_2:=\{(x,\lambda_1,\lambda_2)\in\R^n\times\R\times\R\;|\;\lambda_2\ge f_2(x)\}.
\end{align*}
We can easily verify by \eqref{by def sub} and the normal cone definition that
$$
(v,-1,-1)\in N((\bar{x},f_1(\bar{x}),f_2(\bar{x}));\Omega_1\cap\Omega_2).
$$
To apply Theorem~\ref{nc} to these sets, let us check that $\mbox{\rm ri}(\O_1)\cap\mbox{\rm ri}(\O_2)\ne\emp$. Indeed, we get
\begin{align*}
\mbox{\rm ri}(\Omega_1)&=\{(x,\lambda_1,\lambda_2)\in\R^n\times\R\times\R\;|\; x\in\ri(\dom(f_1)),\;\lambda_1>f_1(x)\}=\ri(\epi(f_1) \times\R,\\
\mbox{\rm ri}(\Omega_2)&=\{(x,\lambda_1,\lambda_2)\in\R^n\times\R\times\R\;|\;x\in\ri(\dom(f_2)),\;\lambda_2>f_2(x)\}
\end{align*}
by Proposition~\ref{re epi}. Then choosing $z\in\ri(\dom(f_1))\cap\ri(\dom(f_2))$, it is not hard to see that
$$
(z,f_1(z)+1,f_2(z)+1)\in\mbox{\rm ri}(\O_1)\cap\mbox{\rm ri}(\O_2)\ne\emp.
$$
Applying now Theorem~\ref{nc} to the above set intersection gives us
$$
N((\bar{x},f_1(\bar{x}),f_2(\bar{x}));\Omega_1\cap\Omega_2)=N((\bar{x},f_1(\bar{x}),f_2(\bar{x}));\Omega_1)+N((\bar{x},f_1(\bar{x}),f_2(\bar{x}));\Omega_2).
$$
It follows from the structures of the sets $\O_1$ and $\O_2$ that
$$
(v,-1,-1)=(v_1,-\gamma_1,0)+(v_2,0,-\gamma_2)
$$
with $(v_1,-\gamma_1)\in N((\bar{x},f_1(\bar{x}));\epi(f_1))$ and $(v_2,-\gamma_2)\in N((\bar{x},f_2(\bar{x}));\epi(f_2))$. Thus
$$
v=v_1+v_2,\quad\gamma_1=\gamma_2=1,
$$
and we have by Proposition~\ref{snr}(ii) that $v_1\in\partial f_1(\bar{x})$ and $v_2\in\partial f_2(\bar{x})$. This ensures the inclusion $\partial(f_1+f_2)(\bar{x})\subset\partial f_1(\bar{x})+\partial f_2(\bar{x})$ and hence verifies (\ref{ssr}) in the case of $m=2$. To complete the proof of the theorem in the general case of $m>2$, we proceed by induction similarly to the proof of Theorem~\ref{nc} with the usage of Lemma~\ref{ri} to deal with relative interiors in the qualification condition (\ref{riq}). $\h$

The next result is a simple consequence of Theorem~\ref{sr} providing a mild condition for the existence of subgradients of an extended-real-valued convex function at a given point.

\begin{Corollary}\label{nt} Let $f\colon\R^n\to(-\infty,\infty]$ be a convex function. Then the validity of the relative interiority condition $\ox\in\mbox{\rm ri(dom}\,f)$ ensures that $\partial f(\ox)\ne\emp$.
\end{Corollary}
{\bf Proof.} Define the extended-real-valued function on $\R^n$ by
\begin{equation*}
g(x):=f(x)+\delta_{\{\ox\}}(x)=\begin{cases}
f(\ox)&\text{if }\;x=\ox,\\
\infty&\text{otherwise}
\end{cases}
\end{equation*}
via the indicator function of the singleton $\{\ox\}$. Then $\epi(g)=\{\ox\}\times[f(\ox),\infty)$ and hence $N((\ox,g(\ox));\epi(g))=\R^n\times(-\infty, 0]$. We obviously get that $\partial g(\ox)=\R^n$ and that $\partial\delta_{\{\ox\}}(\ox)=N(\ox;\{\ox\})=\R^n$ by Proposition~\ref{snr}(i). We further have
\begin{equation*}
\mbox{\rm ri(dom}(h))=\{\ox\}\;\mbox{ for }\;h(x):=\delta_{\{\ox\}}(x)
\end{equation*}
and thus $\mbox{\rm ri(dom}(f))\cap\mbox{\rm ri(dom}(h))\ne\emp$. Applying the subdifferential sum rule of Theorem~\ref{sr} to the above function $g(x)$ at $\ox$ gives us
\begin{equation*}
\R^n=\partial g(\ox)=\partial f(\ox)+\R^n,
\end{equation*}
which justifies the claimed assertion on $\partial f(\ox)\ne\emp$. $\h$

\section{Subdifferential Chain Rule}
\setcounter{equation}{0}

In this section we employ the intersection rule of Theorem~\ref{nc} to derive a chain rule of the subdifferential of a composition of an extended-real-valued function and an affine mapping under which we obviously keep convexity. First we present the following useful lemma.

\begin{Lemma}\label{normal cone to graph} Let $B:\R^n\to\R^p$ be an affine mapping given by $B(x):=Ax+b$, where $A$ is a $p\times n$ matrix and $b\in\R^p$. Then for any $(\ox,\oy)\in{\rm gph}(B)$ we have
\begin{equation*}
N\big((\ox,\oy);\gph(B)\big)=\big\{(u,v)\in\R^n\times\R^p\;\big|\;u=-A^\top v\big\}.
\end{equation*}
\end{Lemma}
{\bf Proof.} It is clear that $\gph(B)$ is convex and $(u,v)\in N((\ox,\oy)\;\gph(B))$ if and only if
\begin{equation}\label{mr1}
\la u,x-\ox\ra+\la v,B(x)-B(\ox)\ra\le 0\;\mbox{ for all }\;x\in\R^n.
\end{equation}
It follows directly from the definitions that
\begin{align*}
&\la u,x-\ox\ra+\la v,B(x)-B(\ox)\ra=\la u,x-\ox\ra+\la v,A(x)-A(\ox)\ra\\
&=\la u,x-\ox\ra +\la A^\top v,x-\ox\ra=\la u+A^\top v,x-\ox\ra.
\end{align*}
This implies the equivalence of (\ref{mr1}) to $\la u+A^\top v,x-\ox\ra\le 0$ for all $x\in \R^n$, and so to $u=-A^\top v$. $\h$

\begin{Theorem}\label{chain rule} Let $f\colon\R^p\to(-\infty, \infty]$ be a convex function, and let $B:\R^n\to\R^p$ be as in Lemma~{\rm\ref{normal cone to graph}} with $B(\ox)\in\dom(f)$ for some $\ox\in\R^n$. Denote $\oy:=B(\ox)$ and assume that the range of $B$ contains a point of $\mbox{\rm ri(dom}(f))$.
Then we have the subdifferential chain rule
\begin{equation}\label{chain-r1}
\partial(f\circ B)(\ox)=A^\top\big(\partial f(\oy)\big)=\big\{A^\top v\;\big|\;v\in\partial f(\oy)\big\}.
\end{equation}
\end{Theorem}
{\bf Proof.} Fix $v\in\partial(f\circ B)(\ox)$ and form the subsets of $\R^n\times\R^p\times\R$ by
\begin{equation*}
\O_1:=\gph(B)\times\R\;\mbox{ and }\;\O_2:=\R^n\times\epi(f).
\end{equation*}
Then we clearly get the relationships
\begin{align*}
\mbox{\rm ri}(\O_1)=\O_1=\gph(B)\times\R,\quad\mbox{\rm ri}(\O_2)=\{(x,y,\lambda)\;|\;x\in\R^n,\;y\in\mbox{\rm ri(dom}(f)),\;\lambda>f(y)\},
\end{align*}
and thus the assumption of the theorem tells us that $\mbox{\rm ri}(\O_1)\cap \mbox{\rm ri}(\O_2)\neq\emp$.

Further, it follows from the definitions of the subdifferential and of the normal cone that $(v,0,-1)\in N((\ox,\oy,\oz);\O_1\cap\O_2)$, where $\oz:=f(\oy)$. Indeed, for any $(x,y,\lambda)\in\O_1\cap\O_2$ we have $y=B(x)$ and $\lambda\ge f(y)$, and so $\lambda\ge f(B(x))$. Thus
\begin{equation*}
\la v,x-\ox\ra+0(y-\oy)+(-1)(\lambda-\oz)\le\la v, x-\ox\ra-[f(B(x))-f(B(\ox))]\le 0.
\end{equation*}
Employing the intersection rule of Theorem~\ref{nc} to the above sets gives us
\begin{equation*}
(v,0,-1)\in N\big((\ox,\oy,\oz);\O_1\big)+N\big((\ox,\oy,\oz);\O_2\big),
\end{equation*}
which reads that $(v,0,-1)=(v,-w,0)+(0,w,-1)$ with $(v,-w)\in N((\ox,\oy);\gph(B))$ and $(w,-1)\in N((\oy,\oz);\epi(f))$. Then we get
\begin{equation*}
v=A^\top w\;\mbox{ and }\;w\in\partial f(\oy),
\end{equation*}
which implies in turn that $v\in A^\top(\partial f(\oy))$ and hence verifies the inclusion ``$\subset$" in \eqref{chain-r1}. The opposite inclusion follows directly from the definition of the subdifferential. $\h$

\section{Subdifferentiation of Maximum Functions}
\setcounter{equation}{0}

Our next topic is subdifferentiation of an important class of nonsmooth convex functions defined as the pointwise maximum of convex functions. We calculate the subdifferential of such functions by using again the geometric intersection rule of Theorem~\ref{nc}.

Given $f_i\colon\R^n\to(-\infty,\infty]$ for $i=1,\ldots,m$, define the {\em maximum function} by
\begin{equation}\label{max-fun}
f(x):=\disp\max_{i=1,\ldots,m}f_i(x),\quad x\in\R^n,
\end{equation}
and for $\ox\in\R^n$ consider the \emph{active index} set
\begin{equation*}
I(\ox):=\big\{i\in\{1,\ldots,m\}\;\big|\;f_i(\ox)=f(\ox)\big\}.
\end{equation*}

\begin{Lemma}\label{ls} {\rm\bf (i)} Let $\O$ be a convex set in $\R^n$. Then $\mbox{\rm int}\,(\O)=\ri(\O)$ provided that $\mbox{\rm int}\,(\O)\ne\emp$. Furthermore, $N(\ox;\O)=\{0\}$ if $\ox\in\mbox{\rm int}(\O)$.\\
{\rm\bf (ii)} Let $f\colon\R^n\to(-\infty,\infty]$ be a convex function, which is continuous at $\ox\in\dom(f)$. Then we have $\ox\in\mbox{\rm int}(\dom(f))$ with the implication
$$
(v,-\lambda)\in N((\ox,f(\ox));\epi(f))\Longrightarrow[\lambda\ge 0\;\mbox{ and }\;v\in\lambda\partial f(\ox)].
$$
\end{Lemma}
{\bf Proof.} {\bf (i)} Suppose that $\mbox{\rm int}(\O)\neq\emptyset$ and check that $\aff(\O)=\R^n$. Indeed, picking $\ox\in\mbox{\rm int}(\O)$ and fixing $x\in \R^n$, find $t>0$ with $tx+(1-t)\ox=\ox+t(x-\ox)\in\mbox{\rm int}(\O)\subset\aff(\O)$. It yields
\begin{equation*}
x=\frac{1}{t}(tx+(1-t)\ox)+(1-\frac{1}{t})\ox\in\mbox{\rm aff}(\O),
\end{equation*}
which justifies the claimed statement due to the definition of relative interior.

To verify the second statement in (i), take $v\in N(\ox;\O)$ with $\ox\in\mbox{\rm int}(\O)$ and get
$$
\la v,x-\ox\ra\le 0\;\mbox{\rm for all }\;x\in\O.
$$
Choosing $\delta>0$ such that $\ox+tv\in\B(\ox;\delta)\subset\O$ for $t>0$ sufficiently small, gives us
\begin{equation*}
\la v, \ox+tv-\ox\ra=t\|v\|^2 \leq 0,
\end{equation*}
which implies $v=0$ and thus completes the proof of assertion (i).\\[1ex]
{\bf (ii)} The continuity of $f$ allows us to find $\delta>0$ such that
$$
|f(x)-f(\ox)|<1\;\mbox{\rm whenever }\;x\in\B(\ox;\delta).
$$
This yields $\B(\ox;\delta)\subset\dom(f)$ and shows therefore that $\ox\in\mbox{\rm int}(\dom(f))$.

Now suppose that $(v,-\lambda)\in N((\ox,f(\ox));\epi(f))$. Then
\begin{equation}\label{nre}
\la v,x-\ox\ra-\lambda(t-f(\ox))\le 0\;\mbox{\rm whenever }\;(x,t)\in\epi(f).
\end{equation}
Employing this inequality with $x=\ox$ and $t=f(\ox)+1$ yields $\lambda\ge 0$.

If $\lambda>0$, we readily get $(v/\lambda,-1)\in N((\ox,f(\ox));\epi(f))$. It follows from Proposition~\ref{snr} that $v/\lambda\in\partial f(\ox)$, and hence $v\in\lambda \partial f(\ox)$.

In the case where $\lambda=0$, we deduce from \eqref{nre} that $v\in N(\ox;(\dom(f))=\{0\}$, and so the inclusion $v\in \lambda\partial f(\ox)$ is also valid. Note that $\partial f(\ox)\ne\emp$ by Corollary~\ref{nt}. $\h$\vspace*{0.05in}

Now we are ready to derive the following {\em maximum rule}.

\begin{Theorem}\label{maxcor} Let $f_i\colon\R^n\to(-\infty,\infty]$, $i=1,\ldots,m$, be convex functions, and let $\ox\in\bigcap_{i=1}^m\dom f_i$ be such that each $f_i$ is continuous at $\ox$. Then we have the maximum rule:
\begin{equation*}
\partial\big(\max f_i)(\ox)={\rm co}\left(\bigcup_{i\in I(\ox)}\partial f_i(\ox)\right).
\end{equation*}
\end{Theorem}
{\bf Proof.} Let $f$ be the maximum function defined in \eqref{max-fun} for which we obviously have
\begin{equation*}
\epi(f)=\bigcap_{i=1}^m\epi(f_i).
\end{equation*}
Employing Proposition~\ref{re epi} and Lemma~\ref{ls}(i) give us the equalities
\begin{align*}
\ri(\epi(f_i))=\{(x,\lambda)\;|\;x\in\mbox{\rm ri}(\dom(f_i)),\;\lambda>f_i(x)\}=\{(x,\lambda)\;|\;x\in\mbox{\rm int}(\dom(f_i)),\;\lambda>f_i(x)\},
\end{align*}
which imply that $(\ox, f(\ox)+1)\in \bigcap_{i=1}^m\mbox{\rm int}(\epi(f_i))=\bigcap_{i=1}^m\ri(\epi(f_i))$. Furthermore, since $f_i(\ox)<f(\ox)=\bar\alpha$ for any $i\notin I(\ox)$, there exists a neighborhood $U$ of $\ox$ and $\gamma>0$ such that $f_i(x)<\alpha$ whenever $(x,\alpha)\in U\times(\bar\alpha-\gamma,\bar\alpha+\gamma)$. It follows that $(\ox,\bar\alpha)\in{\rm int}(\epi(f_i))$, and so $N((\ox, \bar\alpha);\epi(f_i))=\{(0,0)\}$ for such indices $i$. Thus Theorem~\ref{nc} tells us that
\begin{equation*}
N\big((\ox,f(\ox));\epi(f)\big)=\sum_{i=1}^m N\big((\ox,\bar\alpha);\epi(f_i)\big)=\sum_{i\in I(\ox)}N\big((\ox,f_i(\ox));\epi(f_i)\big).
\end{equation*}

Picking now $v\in\partial f(\ox)$, we have by Proposition~\ref{snr}(ii) that $(v,-1)\in N((\ox,f(\ox));\epi f)$, which allows us to find  $(v_i,-\lambda_i)\in N((\ox,f_i(\ox));\epi f_i)$ for $i\in I(\ox)$ such that
\begin{equation*}
(v,-1)=\sum_{i\in I(\ox)}(v_i,-\lambda_i).
\end{equation*}
This yields $\sum_{i\in I(\ox)}\lambda_i=1$, $\lambda_i\ge 0$, $v=\sum_{i\in I(\ox)}v_i$, and $v_i\in\lambda_i\partial f_i(\ox)$ by Lemma \ref{ls}(ii). Thus
$v=\sum_{i\in I(\ox)}\lambda_i u_i$, where $u_i\in\partial f_i(\ox)$ and $\sum_{i\in I(\ox)}\lambda_i=1$. This verifies that
\begin{equation*}
v\in {\rm co}\left(\bigcup_{i\in I(\ox)}\partial f_i(\ox)\right).
\end{equation*}
The opposite inclusion in the maximum rule follows from
$$
\partial f_i(\ox)\subset\partial f(\ox)\;\mbox{for all }\;i\in I(\ox),
$$
which in turn follows directly from the definitions. $\h$

\section{Optimal Value Function and Another Chain Rule}
\setcounter{equation}{0}

The main result of this section concerns calculating the subdifferential of extended-real-valued convex functions, which play a remarkable role in variational analysis, optimization, and their numerous applications and are known under the name of {\em optimal value/marginal functions}. Functions of this class are generally defined by
\begin{equation}\label{optimal value}
\mu(x):=\inf\big\{\ph(x,y)\;\big|\;y\in F(x)\big\},
\end{equation}
where $\ph\colon\R^n\times\R^p\to(-\infty,\infty]$ is an extended-real-valued function, and where $F\colon\R^n\tto\R^p$ is a set-valued mapping, i.e., $F(x)\subset\R^p$ for every $x\in\R^n$. In what follows we select $\ph$ and $F$ in such a way that the resulting function (\ref{optimal value}) is convex and to derive a formula to express its subdifferential via the subdifferential of $\ph$ and an appropriate generalized differentiation construction for the set-valued mapping $F$. The results obtained in the general framework of variational analysis \cite{m-book1,rw} advise us that the most suitable construction for $F$ for these purposes is the so-called {\em coderivative} of $F$ at $(\ox,\oy)\in\gph(F)$ defined via the normal cone to the graphical set $\gph(F):=\{(x,y)\in\R^n\times\R^p\;|\;y\in F(x)\}$ by
\begin{equation}\label{cod}
D^*F(\ox,\oy)(v)=\{u\in\R^n\;|\;(u,-v)\in N((\ox,\oy);\gph(F))\},\quad v\in\R^p.
\end{equation}
It is easy to check that the optimal value function (\ref{optimal value}) is convex provided that $\ph$ is convex and the graph of $F$ is convex as well. An example of such $F$ is given by the affine mapping $B$ considered in Lemma~\ref{normal cone to graph}. Note that, as follows directly from Lemma~\ref{normal cone to graph} and definition (\ref{cod}), the coderivative of this mapping is calculated by
\begin{eqnarray}\label{cod1}
D^*B(\ox,\oy)(v)=A^\top v\;\mbox{\rm with }\;\oy=B(\ox).
\end{eqnarray}

Now we proceed with calculating the subdifferential of (\ref{optimal value}) via that of $\ph$ and the coderivative of $F$. The results of this type are induced by those in variational analysis \cite{m-book1,rw}, where only upper estimate of $\partial\mu(\ox)$ were obtained. The convexity setting of this paper and the developed approach allow us to derive an {\em exact formula} (equality) for calculating $\partial\mu(\ox)$ under a mild relative interior condition, which is strictly weaker than the normal qualification condition from \cite[Theorem~2.61]{bmn}; cf.\ the discussion at the end of Section~5.

\begin{Theorem}\label{mr3} Let $\mu(\cdot)$ be the optimal value function {\rm(\ref{optimal value})} generated by a convex-graph mapping $F:\R^n\tto\R^p$ and a convex function $\ph:\R^n\times\R^p\to(-\infty,\infty]$. Suppose that $\mu(x)>-\infty$ for all $x\in\R^n$, fix some $\ox\in\dom(\mu)$, and consider the solution set
\begin{equation*}
S(\ox):=\big\{\oy\in F(\ox)\;\big|\;\mu(\ox)=\ph(\ox,\oy)\big\}.
\end{equation*}
If $S(\ox)\ne\emp$, then for any $\oy\in S(\ox)$ we have the equality
\begin{equation}\label{vf2}
\partial\mu(\ox)=\bigcup_{(u,v)\in\partial\ph(\ox,\oy)}\big[u+D^*F(\ox,\oy)(v)\big]
\end{equation}
provided the validity of the relative interior qualification condition
\begin{equation}\label{qf}
\mbox{\rm ri(dom}(\ph))\cap\mbox{\rm ri(\gph}(F))\ne\emp.
\end{equation}
\end{Theorem}
{\bf Proof.} Picking any $\oy\in S(\ox)$, let us first verify the estimate
\begin{equation}\label{vf1}
\bigcup_{(u,v)\in\partial\ph(\ox,\oy)}\big[u+D^*F(\ox,\oy)(v)\big]\subset\partial\mu(\ox).
\end{equation}
To proceed, take $w$ from the set on the left-hand side of \eqref{vf1} and find $(u,v)\in\partial\ph(\ox,\oy)$ with $w-u\in D^*F(\ox,\oy)(v)$.
It gives us $(w-u,-v)\in N((\ox,\oy);\gph(F))$ and thus
\begin{eqnarray*}
\la w-u,x-\ox\ra-\la v,y-\oy\ra\le 0\;\mbox{ for all }\;(x,y)\in\gph(F),
\end{eqnarray*}
which shows that whenever $y\in F(x)$ we have
\begin{equation*}
\la w,x-\ox\ra\le\la u,x-\ox\ra+\la v,y-\oy\ra\le\ph(x,y)-\ph(\ox,\oy)=\ph(x,y)-\mu(\ox).
\end{equation*}
This allows us to arrive at the estimate
\begin{equation*}
\la w,x-\ox\ra\le\inf_{y\in F(x)}\ph(x,y)-\mu(\ox)=\mu(x)-\mu(\ox)
\end{equation*}
justifying the inclusion $w\in\partial\mu(\ox)$ and hence the claimed one in \eqref{vf1}.

It remains to verify the inclusion ``$\subset$" in \eqref{vf2}. Take $w\in\partial\mu(\ox)$, $\oy\in S(\ox)$ and get
\begin{align*}
\la w,x-\ox\ra&\le\mu(x)-\mu(\ox)=\mu(x)-\ph(\ox,\oy)\le\ph(x,y)-\ph(\ox,\oy)
\end{align*}
whenever $y\in F(x)$ and $x\in\R^n$. This implies in turn that for any $(x,y)\in\R^n\times\R^p$ we have
\begin{equation*}
\la w,x-\ox\ra +\la 0,y-\oy\ra\le\ph(x,y)+\delta\big((x,y);\gph(F))-\big[\ph(\ox,\oy)+\delta\big((\ox,\oy);\gph(F)\big)\big].
\end{equation*}
Considering further $f(x,y):=\ph(x,y)+\delta((x,y);\gph(F))$, deduce from the subdifferential sum rule of Theorem~\ref{sr} under \eqref{qf} that
\begin{equation*}
(w,0)\in\partial f(\ox,\oy)=\partial\ph(\ox,\oy)+N\big((\ox,\oy);\gph(F)\big).
\end{equation*}
This shows that $(w,0)=(u_1,v_1)+(u_2,v_2)$ with $(u_1,v_1)\in\partial\ph(\ox,\oy)$ and $(u_2,v_2)\in N((\ox,\oy);\gph(F))$ and thus yields $v_2=-v_1$. Hence $(u_2,-v_1)\in N((\ox,\oy);\gph(F))$ meaning by definition that $u_2\in D^*F(\ox,\oy)(v_1)$. Therefore we arrive at
\begin{equation*}
w=u_1+u_2\in u_1+D^*F(\ox,\oy)(v_1),
\end{equation*}
which justifies the inclusion ``$\subset$" in (\ref{vf2}) and completes the proof of the theorem. $\h$

Observe that Theorem~\ref{mr3} easily implies the chain rule of Theorem~\ref{chain rule} by setting $F(x):=\{B(x)\}$ and $\ph(x,y):=f(y)$ therein. Then we have  $\mu(x)=(f\circ B)(x)$,
\begin{equation*}
\mbox{\rm ri(dom}(\ph))=\R^n\times\mbox{\rm ri(dom}(f)),\;\mbox{\rm ri(\gph}(F))=\gph(B),
\end{equation*}
Thus the relative interiority assumption  of Theorem~\ref{chain rule} yields the validity of the qualification condition (\ref{qf}) is satisfied, and we arrive at the chain rule (\ref{chain-r1}) directly from (\ref{vf2}) and the coderivative expression in (\ref{cod1}).

We now derive from Theorem~\ref{mr3} and the intersection rule of Theorem~\ref{nc} a new subdifferential chain rule concerning compositions of convex functions with particular structures. We say that $g:\R^p\to(-\infty,\infty]$ is \emph{nondecreasing componentwise} if
\begin{equation*}
\big[x_i\le y_i\;\mbox{ for all }\;i=1,\ldots,p\big]\Longrightarrow\big[g(x_1,\ldots,x_p)\le g(y_1,\ldots,y_p)\big].
\end{equation*}

\begin{Theorem}\label{general com} Define $h:\R^n\to\R^p$ by $h(x):=(f_1(x),\ldots,f_p(x))$, where $f_i:\R^n\to\R$ for $i=1,\ldots,p$ are convex functions. Suppose that $g:\R^p\to(-\infty,\infty]$ is convex and nondecreasing componentwise. Then the composition $g\circ h:\R^n\to(-\infty,\infty]$ is a convex function, and we have the subdifferential chain rule
\begin{equation}\label{cwc}
\partial(g\circ h)(\ox)=\Big\{\sum_{i=1}^p\gamma_i v_i\;\Big|\;(\gamma_1,\ldots,\gamma_p)\in\partial g(\oy),\;v_i\in\partial f_i(\ox),\;i=1,\ldots,p\Big\}
\end{equation}
with $\ox\in\R^n$ and $\oy:=h(\ox)\in \dom(g)$ under the condition that there exist $\ou\in\R^n$ and $\bar\lambda_i>f_i(\ou)$ for all $i=1,\ldots,p$ satisfying
$$
(\bar\lambda_1,\ldots,\bar\lambda_p)\in\ri(\dom (g)).
$$
\end{Theorem}
{\bf Proof.} Let $F:\R^n\tto\R^p$ be a set-valued mapping defined by
\begin{equation*}
F(x):=[f_1(x),\infty)\times[f_2(x),\infty)\times\ldots\times[f_p(x),\infty).
\end{equation*}
Then the graph of $F$ is represented by
\begin{equation*}
\gph(F)=\big\{(x,t_1,\ldots,t_p)\in\R^n\times\R^p\big|\;t_i\ge f_i(x)\big\}.
\end{equation*}
Consider further the convex sets
\begin{equation*}
\O_i:=\{(x,\lambda_1,\ldots,\lambda_p)\;|\;\lambda_i\ge f_i(x)\}
\end{equation*}
and observe that $\gph(F)=\bigcap_{i=1}^p\O_i$.
Since all the functions $f_i$ are convex, the set $\gph(F)$ is convex as well. Define $\ph:\R^n\times\R^p\to(-\infty,\infty]$ by $\ph(x,y):=g(y)$ and observe, since $g$ is increasing componentwise, that
\begin{equation*}
\inf\big\{\ph(x,y)\;\big|\;y\in F(x)\big\}=g\big(f_1(x),\ldots,f_p(x)\big)=(g\circ h)(x),
\end{equation*}
which ensures the convexity of the composition $g\circ h$; see \cite[Proposition~1.54]{bmn}. It follows from Proposition~\ref{re epi} and Lemma~\ref{ri} that
\begin{equation*}
\ri(\gph(F))=\{(x,\lambda_1,\ldots,\lambda_p)\;|\;\lambda_i>f_i(x)\;\mbox{\rm for all }\;i=1,\ldots,p\}.
\end{equation*}
The assumptions made in the theorem guarantee that $\ri(\gph(F))\cap\ri(\dom(\ph))\ne\emp$.  Moreover, the structure of each set $\O_i$ gives us
$$
[(v,-\gamma_1,\ldots,-\gamma_p)\in N((\ox,f_1(\ox),\ldots,f_p(\ox));\O_i)]\Longleftrightarrow[(v,-\gamma_i)\in N((\ox, f_i(\ox));\epi(f_i)),\gamma_j=0\;\mbox{if } j\ne i].
$$
Using the coderivative definition and applying the intersection rule of Theorem~\ref{nc} and also Lemma~\ref{ls}, we get
\begin{align*}
v\in D^*F(\ox, f_1(\ox),\ldots,f_p(\ox))(\gamma_1,\ldots,\gamma_p)&\Longleftrightarrow(v,-\gamma_1,\ldots,-\gamma_p)\in N((\ox, f_1(\ox),\ldots,f_p(\ox));\gph(F)))\\
&\Longleftrightarrow(v,-\gamma_1,\ldots,-\gamma_p)\in N((\ox,f_1(\ox),\ldots,f_p(\ox));\bigcap_{i=1}^p\O_i)\\
&\Longleftrightarrow(v,-\gamma_1,\ldots,-\gamma_p)\in\sum_{i=1}^p N((\ox,f_1(\ox),\ldots,f_p(\ox));\O_i)\\
&\Longleftrightarrow v=\sum_{i=1}^pv_i\;\mbox{\rm with }\;(v_i,-\gamma_i)\in N((\ox, f_i(\ox));\epi f_i)\\
&\Longleftrightarrow v=\sum_{i=1}^pv_i\;\mbox{\rm with }\;v_i\in\gamma_i\partial f_i(\ox)\\
&\Longleftrightarrow v\in\sum_{i=1}^p\gamma_i\partial f_i(\ox).
\end{align*}
It follows from Theorem~\ref{mr3} that $v\in\partial(g\circ h)(\ox)$ if and only if there exists a collection $(\gamma_1,\ldots,\gamma_p)\in\partial g(\oy)$ such that $v\in D^*F(\ox,f_1(\ox),\ldots,f_p(\ox))(\gamma_1,\ldots,\gamma_p)$. This allows us to deduce the chain rule \eqref{cwc} from the equivalences above. $\h$

\section{Normals to Preimages of Sets via Set-Valued Mappings}
\setcounter{equation}{0}

In this section we calculate the normal cone to convex sets of a special structure that frequently appear in variational analysis and optimization. These sets are constructed as follows. Given a set $\Theta\subset\R^p$ and a set-valued mapping $F\colon\R^n\tto\R^p$ the {\em preimage} or {\em inverse image} of the set $\Th$ under the mapping $F$ is defined by
\begin{equation}\label{pre}
F^{-1}(\Theta):=\{x\in\R^n\;|\;F(x)\cap\Theta\ne\emp\}.
\end{equation}
Our goal here is to calculate the normal cone to the preimage set (\ref{pre}) via the normal cone to $\Th$ and the coderivative of $F$. This is done in the following theorem, which is yet another consequence of the intersection rule from Theorem~\ref{nc}.

\begin{Theorem} Let $F\colon\R^n\tto\R^p$ be a set-valued mapping with convex graph, and let $\Theta$ be a convex subset of $\R^p$. Suppose that there exists $(a,b)\in\R^n\times\R^p$ satisfying
\begin{equation*}
(a,b)\in\ri(\gph(F))\;\mbox{\rm and }\;b\in\ri(\Theta).
\end{equation*}
Then for any $\ox\in F^{-1}(\Theta)$ and $\oy\in F(\ox)\cap\Theta$ we have the representation
\begin{equation}\label{nci}
N(\ox;F^{-1}(\Theta))=D^*F(\ox,\oy)(N(\oy;\Theta)).
\end{equation}
\end{Theorem}
{\bf Proof.} It is not hard to show that $F^{-1}(\Theta)$ is a convex set. Picking any $u\in N(\ox;F^{-1}(\Theta))$ gives us by definition that
\begin{equation*}
\la u,x-\ox\ra\le 0\;\mbox{\rm whenever }\;x\in F^{-1}(\Theta),\;\mbox{\rm i.e., }\;F(x)\cap\Theta\ne\emp.
\end{equation*}
Consider the two convex subsets of $\R^{n+p}$ defined by
\begin{equation*}
\O_1:=\gph (F)\;\mbox{\rm and }\;\O_2:=\R^n\times\Theta
\end{equation*}
for which we have $(u,0)\in N((\ox,\oy);\O_1\cap\O_2)$. Applying now Theorem~\ref{nc} tells us that
\begin{equation*}
(u,0)\in N((\ox,\oy);\O_1)+N((\ox,\oy);\O_2)=N((\ox,\oy);\gph F)+[\{0\}\times N(\oy;\Theta)],
\end{equation*}
and thus we get the representation
\begin{equation*}
(u,0)=(u,-v)+(0,v)\;\mbox{ with }\;(u,-v)\in N((\ox,\oy);\gph F)\;\mbox{ and }\;v\in N(\oy;\Theta)
\end{equation*}
from which it follows immediately that
\begin{equation*}
u\in D^*F(\ox,\oy)(v)\;\mbox{\rm and }\;v\in N(\oy;\Theta).
\end{equation*}
This versifies the inclusion``$\subset$" in (\ref{nci}). The opposite inclusion is trivial. $\h$

\section{Coderivative Calculus}
\setcounter{equation}{0}

We see from above that the coderivative notion (\ref{cod}) is instrumental to deal with set-valued mappings. Although this notion was not properly developed in basic convex analysis, the importance of it has been fully revealed in general variational analysis and its applications; see, e.g., \cite{bz,m-book1,rw} and the references therein, where the reader can find, in particular, various results on coderivative calculus. Most of these results were obtained in the inclusion form under the corresponding normal qualification conditions generated by (\ref{bqc}). We present below some calculus rules for coderivatives of convex-graph mappings, which are derived from the intersection rule of Theorem~\ref{nc} and hold as {\em equalities} in a bit different form under weaker {\em relative interior qualification conditions}.

Recall that the domain of a set-valued mapping $F\colon\R^n\tto\R^p$ is defined by
$$
\dom(F):=\{x\in\R^m\;|\;F(x)\ne\emp\}.
$$
Given two set-valued mappings $F_1,F_2\colon\R^n\tto\R^p$, their sum is defined by
\begin{equation*}
(F_1+F_2)(x)=F_1(x)+F_2(x):=\{y_1+y_2\;|\;y_1\in F_1(x),\;y_2\in F_2(x)\}.
\end{equation*}
It is easy to see that $\dom(F_1+F_2)=\dom(F_1)\cap\dom(F_2)$ and that the graph of $F_1+F_2$ is convex provided that both $F_1,F_2$ have this property.

Our first calculus result concerns representing the coderivative of the sum $F_1+F_2$ at the given point $(\ox,\oy)\in\gph(F_1+F_2)$. To formulate it, consider the nonempty set
\begin{equation*}
S(\ox,\oy):=\left\{(\oy_1,\oy_2)\in\R^p\times\R^p\;|\;\oy=\oy_1+\oy_2,\;\oy_i\in F_i(\ox)\;\mbox{\rm for }\;i=1,2\right\}.
\end{equation*}

\begin{Theorem} Let $F_1,F_2\colon\R^n\tto\R^p$ be set-valued mappings of convex graphs, and let the relative interior qualification condition
\begin{equation}\label{QC}
\ri(\gph(F_1))\cap\ri(\gph(F_2))\ne\emp
\end{equation}
hold. Then we have the coderivative sum rule
\begin{equation}\label{csr}
D^*(F_1+F_2)(\ox,\oy)(v)=\bigcap_{(\oy_1,\oy_2)\in S(\ox,\oy)}\left[D^*F_1(\ox,\oy_1)(v)+D^*F_2(\ox,\oy_2)(v)\right]
\end{equation}
for all $(\ox,\oy)\in\gph(F_1+F_2)$ and $v\in\R^p$.
\end{Theorem}
{\bf Proof.} Fix any $u\in D^*(F_1+F_2)(\ox,\oy)(v)$ and $(\oy_1,\oy_2)\in S(\ox,\oy)$ for which we have the inclusion $(u,-v)\in N((\ox,\oy);\gph(F_1+F_2))$.
Consider the convex sets
\begin{align*}
\O_1:=\{(x,y_1,y_2)\in \R^n\times\R^p\times\R^p\;|\;y_1\in F_1(x)\},\quad\O_2:=\{(x,y_1,y_2)\in\R^n\times\R^p\times\R^p\;|\;y_2\in F_2(x)\}
\end{align*}
and deduce from the normal cone definition that
\begin{equation*}
(u,-v,-v)\in N((\ox,\oy_1,\oy_2);\O_1\cap\O_2).
\end{equation*}
It is easy to observe the relative interior representations
\begin{align*}
&\ri(\O_1)=\{(x,y_1,y_2)\in\R^n\times\R^p\times\R^p\;|\;(x,y_1)\in\ri(\gph (F_1))\},\\
&\ri(\O_2)=\{(x,y_1,y_2)\in\R^n\times\R^p\times\R^p\;|\;(x,y_2)\in\ri(\gph (F_2))\},
\end{align*}
which show that condition \eqref{QC} yields $\ri(\O_1)\cap\ri(\O_2)\ne\emp$. It tells us by Theorem~\ref{nc} that
\begin{equation*}
(u,-v,-v)\in N((\ox,\oy_1,\oy_2);\O_1)+N((\ox,\oy_1,\oy_2);\O_2),
\end{equation*}
and thus we arrive at the representation
\begin{equation*}
(u,-v,-v)=(u_1,-v,0)+(u_2,0,-v)\;\mbox{ with }\;(u_i,-v)\in N((\ox,\oy_i);\gph(F_i)),\;i=1,2.
\end{equation*}
This verifies therefore the relationship
$$
u=u_1+u_2\in D^*F_1(\ox,\oy_1)(v)+D^*F_2(\ox,\oy_2)(v),
$$
which justifies the inclusion ``$\subset$" in \eqref{csr}. The opposite inclusion is obvious. $\h$

Next we define the composition of two mappings $F:\R^n\tto\R^p$ and $G:\R^p\tto\R^q$ by
\begin{equation*}
(G\circ F)(x)=\bigcup_{y\in F(x)}G(y):=\{z\in G(y)\;|\;y\in F(x)),\quad x\in\R^n,
\end{equation*}
and observe that $G\circ F$ is convex-graph provided that both $F$ and $G$ have this property. Given $\oz\in(G\circ F)(\ox)$, we consider the set
\begin{equation*}
M(\ox,\oz):=F(\ox)\cap G^{-1}(\oz).
\end{equation*}

\begin{Theorem}\label{scr} Let $F:\R^n\tto\R^p$ and $G:\R^p\tto\R^q$ be set-valued mappings of convex graphs for which there exist vectors $(x,y,z)\in\R^n\times\R^p\times\R^q$ satisfying
\begin{equation}\label{QC1}
(x,y)\in\ri(\gph(F))\;\mbox{\rm and }\;(y,z)\in\mbox{\rm ri}(\gph(G)).
\end{equation}
Then for any $(\ox,\oz)\in\gph(G\circ F)$ and $w\in\R^q$ we have the coderivative chain rule
\begin{equation}\label{chain}
D^*(G\circ F)(\ox,\oz)(w)=\bigcap_{\oy\in M(\ox,\oz)}D^*F(\ox,\oy)\circ D^*G(\oy,\oz)(w).
\end{equation}
\end{Theorem}
{\bf Proof.} Picking $u\in D^*(G\circ F)(\ox,\oz)(w)$ and $\oy\in M(\ox,\oz)$ gives us the inclusion $(u,-w)\in N((\ox,\oz);\gph(G\circ F))$, which means that
\begin{equation*}
\la u,x-\ox\ra-\la w,z-\oz\ra\le 0\;\mbox{\rm for all }\;(x,z)\in\gph(G\circ F).
\end{equation*}
Form now the two convex subsets of $\R^n\times\R^p\times\R^q$ by
$$
\O_1:=\gph(F)\times\R^q\;\mbox{ and }\;\O_2:=\R^n\times\gph(G)
$$
We can easily deduce from the definitions that
\begin{equation*}
(u,0,-w)\in N((\ox,\oy,\oz);\O_1\cap\O_2)
\end{equation*}
and that the qualification condition (\ref{QC1}) ensures the validity of the one $\ri(\O_1)\cap\ri(\O_2)\ne\emp$ in Theorem~\ref{nc}. Applying then the intersection rule to the above sets tells us that
\begin{equation*}
(u,0,-w)\in N((\ox,\oy,\oz);\O_1\cap\O_2)=N((\ox,\oy,\oz);\O_1)+N((\ox,\oy,\oz);\O_2),
\end{equation*}
and thus there is a vector $v\in\R^p$ such that we have the representation
\begin{equation*}
(u,0,-w)=(u, -v,0)+(0,v,-w)\;\mbox{with}\;(u,-v)\in N((\ox,\oy);\gph(F)),\;(v,-w)\in N((\oy,\oz);\gph(G)).
\end{equation*}
This shows by the coderivative definition (\ref{cod}) that
\begin{equation*}
u\in D^*F(\ox,\oy)(v)\;\mbox{\rm and }\;v\in D^*G(\oy,\oz)(w),
\end{equation*}
and so we justify the inclusion ``$\subset$" in \eqref{chain}. The opposite inclusion is easy to verify. $\h$

Our final result in this section provides an exact formula for calculating the coderivative of intersections of set-valued mappings $F_1,F_2:\R^n\tto\R^p$ defined by
\begin{equation*}
(F_1\cap F_2)(x):=F_1(x)\cap F_2(x),\quad x\in\R^n,
\end{equation*}
which is also deduced from the basic intersection rule for the normal cone in Theorem~\ref{nc}.

\begin{Proposition}\label{cointersection} Let $F_1,F_2:\R^n\tto\R^p$ be of convex graphs, and let
\begin{equation*}
\ri(\gph(F_1))\cap\ri(\gph(F_2))\ne\emp.
\end{equation*}
Then for any $\oy\in (F_1\cap F_2)(\ox)$ and $v\in\R^p$ we have
\begin{equation}\label{cim}
D^*(F_1\cap F_2)(v)=\bigcup_{v_1+v_2=v}\left[D^*F_1(\ox,\oy)(v_1)+D^*F_2(\ox,\oy)(v_2)\right].
\end{equation}
\end{Proposition}
{\bf Proof.} It follows from the definition that $\gph(F_1\cap F_2)=\gph(F_1)\cap\gph(F_2)$.
Pick any vector $u\in D^*(F_1\cap F_2)(v)$ and get by Theorem~\ref{nc} that
\begin{equation*}
(u,-v)\in N((\ox,\oy);\gph(F_1\cap F_2))=N((\ox,\oy);\gph(F_1))+N((\ox,\oy);\gph(F_2)).
\end{equation*}
This allows us to represent the pair $(u,v)$ in the form
\begin{equation*}
(u,-v)=(u_1,-v_1)+(u_2,-v_2),
\end{equation*}
where $(u_1,-v_1)\in N((\ox,\oy);\gph(F_1))$ and $(u_2,-v_2)\in N((\ox,\oy);\gph(F_2))$. Therefore we have
\begin{equation*}
u=u_1+u_2\in D^*F_1(\ox,\oy)(v_1)+D^*F_2(\ox,\oy)(v_2)\;\mbox{ with }\;v=v_1+v_2,
\end{equation*}
verifying the inclusion ``$\subset$" in \eqref{cim}. The opposite inclusion comes from the definition. $\h$

\section{Solution Maps for Parameterized Generalized Equations}
\setcounter{equation}{0}

Here we present a rather simple application of coderivative calculus to calculating the coderivative of set-valued mappings given in the structural form
\begin{equation}\label{GES}
S(x)=\{y\in\R^p\;|\;0\in F(x,y)+G(x,y)\},\quad x\in\R^n,
\end{equation}
where $F,G:\R^n\times\R^p\tto\R^q$ are set-valued mappings. Mappings of this type can be treated as solutions maps to the so-called {\em generalized equations}
\begin{equation*}
0\in F(x,y)+G(x,y),\quad x\in\R^n,\;y\in\R^p,
\end{equation*}
with respect to the decision variable $y$ under parameterization/perturbation by $x$. This terminology and first developments go back to Robinson \cite{rob}, who considered the case where $G(y)=N(y;\O)$ is the normal cone mapping associated with a convex set $\O$ and where $F(x,y)$ is single-valued. It has been recognized that the generalized equation formalism, including its extended form, is a convenient model to investigate various aspects of optimization, equilibrium, stability, etc. In particular, the coderivative of the solution map (\ref{GES}) plays an important role in such studies; see, e.g., \cite{m-book1} and the references therein.

To proceed with calculating the coderivative of the solution map (\ref{GES}) in the case of convex-graph set-valued mappings, we first observe the following  fact of its own interest.

\begin{Proposition}\label{domain} Let $F:\R^n\tto\R^p$ be an arbitrary set-valued mapping with convex graph. Given $\ox\in\dom(F)$, we have the relationships
\begin{equation}\label{domain1}
N(\ox;\dom(F))=D^*F(\ox,\oy)(0)\;\mbox{ for every }\;\oy\in F(\ox).
\end{equation}
\end{Proposition}
{\bf Proof.} Picking any $v\in N(\ox;\dom(F))$ and $\oy\in F(\ox)$ gives us
\begin{equation*}
\la v,x-\ox\ra\le 0\;\mbox{\rm for all }\;x\in\dom(F),
\end{equation*}
which immediately implies the inequality
\begin{equation*}
\la v, x-\ox\ra +\la 0,y-\oy\ra\le 0,\quad y\in F(x).
\end{equation*}
This yields in turn $(v,0)\in N((\ox,\oy);\gph)(F))$ and so $v\in D^*F(\ox,\oy)(0)$ thus verifying the inclusion ``$\subset$" in (\ref{domain1}).
The opposite inclusion in (\ref{domain1}) is straightforward. $\h$

\begin{Theorem} Let $F$ and $G$ in {\rm(\ref{GES})} be convex-graph, and let $(\ox,\oy)\in\gph S$. Impose the qualification condition
\begin{equation}\label{geq}
\ri(\gph(F))\cap\ri(-\gph(G))\ne\emp.
\end{equation}
Then for every $\oz\in F(\ox,\oy)\cap[-G(\ox,\oy)]$ we have
\begin{equation}\label{cod-ge}
D^*S(\ox,\oy)(v)=\bigcup_{w\in \R^q}\{u\in\R^n\;|\;(u,-v)\in D^*F((\ox,\oy),\oz)(w)+D^*G((\ox,\oy),-\oz)(w)\}.
\end{equation}
\end{Theorem}
{\bf Proof.} It is easy to see that the solution map $S$ is convex-graph under this property for $F$ and $G$. Furthermore, we get
\begin{align*}
\gph(S)&=\{(x,y)\in \R^n\times\R^p\; |\; 0\in F(x,y)+G(x,y)\}\\
&=\{(x,y)\in \R^n\times\R^p\; |\; F(x,y)\cap [-G(x,y)]\neq\emp\}=\dom(H),
\end{align*}
where $H(x,y):=F(x,y)\cap [-G(x,y)]$. Take further any $u\in D^*S(\ox,\oy)(v)$ and deduce from (\ref{cod}) and Proposition~\ref{domain} that
\begin{equation*}
(u,-v)\in N((\ox,\oy);\gph(S))=N((\ox,\oy);\dom(H))=D^*H((\ox,\oy,\oz))(0)
\end{equation*}
for every $\oz\in H(\ox,\oy)=F(\ox,\oy)\cap[-G(\ox,\oy)]$. Then the coderivative intersection rule of Proposition~\ref{cointersection} tells us that under the validity of (\ref{geq}) that
\begin{align*}
(u,-v)\in D^*H((\ox,\oy,\oz))(0)&=\bigcup_{w\in \R^q} [D^*F((\ox,\oy),\oz)(w)+D^*(-G)((\ox,\oy),\oz)(-w)]\\
&=\bigcup_{w\in\R^q} [D^*F((\ox,\oy),\oz)(w)+D^*(G)((\ox,\oy),-\oz)(w)],
\end{align*}
which yields (\ref{cod-ge}) and thus completes the proof of the theorem. $\h$

{\bf Acknowledgement.} The authors are grateful to two anonymous referees and the handling Editor for their valuable remarks, which allowed us to improve the original presentation.

\end{document}